\documentclass[11pt]{amsart}

\usepackage{amsmath,amsthm}
\usepackage{amssymb,amsfonts}
\usepackage{latexsym}
\usepackage{mathrsfs}
\usepackage[psamsfonts]{eucal}
\usepackage{bm}
\usepackage{times}
\usepackage{mathptmx}
\usepackage[hang]{caption}
\usepackage[T1]{fontenc}
\usepackage{textcomp}

\usepackage{graphicx,epsfig}

\theoremstyle{plain}
\newtheorem{theorem}{Theorem}[section]
\newtheorem{lem}[theorem]{Lemma}
\newtheorem{rem}[theorem]{Remark}

\newtheorem{cor}[theorem]{Corollary}

\theoremstyle{definition}
\newtheorem{definition}[theorem]{Definition}

\theoremstyle{remark}

\newcommand{\af}{almost Fuchsian}
\newcommand{\ab}{asymptotic boundary}
\newcommand{\abs}{asymptotic boundaries}
\newcommand{\cc}{convex core}

\newcommand{\cmc}{constant mean curvature}

\newcommand{\fg}{fundamental group}
\newcommand{\fp}{fundamental polyhedron}

\newcommand{\Hd}{Hausdorff dimension}

\newcommand{\hym}{hyperbolic metric}
\newcommand{\htm}{hyperbolic three-manifold}
\newcommand{\jc}{Jordan curve}
\newcommand{\kg}{Kleinian group}

\newcommand{\mc}{mean curvature}
\newcommand{\mct}{minimal catenoid}

\newcommand{\ms}{minimal surface}

\newcommand{\pc}{principal curvature}
\newcommand{\qf}{quasi-Fuchsian}
\newcommand{\qfg}{quasi-Fuchsian group}
\newcommand{\qfm}{quasi-Fuchsian manifold}
\newcommand{\RS}{Riemann surface}

\newcommand{\TS}{Teichm\"{u}ller space}
\newcommand{\tm}{three-manifold}

\newcommand{\tdc}{two disjoint circles}

\newcommand{\wrt}{with respect to}

\newcommand{\be}{\begin{equation}}
\newcommand{\ene}{\end{equation}}
\newcommand{\br}{\begin{rem}}
\newcommand{\er}{\end{rem}}
\newcommand{\bl}{\begin{lem}}
\newcommand{\el}{\end{lem}}
\newcommand{\bd}{\begin{Def}}
\newcommand{\ed}{\end{Def}}
\newcommand{\ben}{\begin{enumerate}}
\newcommand{\een}{\end{enumerate}}
\newcommand{\bp}{\begin{proof}}
\newcommand{\ep}{\end{proof}}
\newcommand{\bpo}{\begin{pro}}
\newcommand{\epo}{\end{pro}}
\newcommand{\beq}{\begin{equation*}}
\newcommand{\eeq}{\end{equation*}}
\newcommand{\bear}{\begin{eqnarray*}}
\newcommand{\eear}{\end{eqnarray*}}
\newcommand{\bt}{\begin{theorem}}
\newcommand{\et}{\end{theorem}}

\DeclareMathAlphabet{\mathcal}{OMS}{cmsy}{m}{n}

\numberwithin{equation}{section}

\allowdisplaybreaks

\def\XXint#1#2#3{{\setbox0=\hbox{$#1{#2#3}{\int}$}
    \vcenter{\hbox{$#2#3$}}\kern-.5\wd0}}

\makeatletter
\def\@citestyle{\m@th\upshape\mdseries}
\def\citeform#1{{\bfseries#1}}
\def\@cite#1#2{{%
  \@citestyle[\citeform{#1}\if@tempswa, #2\fi]}}
\@ifundefined{cite }{%
  \expandafter\let\csname cite \endcsname\cite
  \edef\cite{\@nx\protect\@xp\@nx\csname cite \endcsname}%
}{}
\makeatother

\renewcommand{\leq}{\leqslant}
\renewcommand{\geq}{\geqslant}

\newcommand{\C}{\mathbb{C}}
\newcommand{\R}{\mathbb{R}}

\newcommand{\B}{\mathbb{B}}
\renewcommand{\H}{\mathbb{H}}
\newcommand{\X}{\mathbf{X}}

\newcommand{\T}{\mathbf{T}}

\DeclareMathOperator{\dist}{dist}%
\DeclareMathOperator{\genus}{genus}%
\DeclareMathOperator{\id}{id}%

\DeclareMathOperator{\Hull}{Hull}%
\DeclareMathOperator{\Int}{Int}%
\DeclareMathOperator{\Isom}{Isom}%
\DeclareMathOperator{\Length}{Length}%
\DeclareMathOperator{\Mob}{M\text{\"o}b}%
\DeclareMathOperator{\PSL}{PSL}%
\DeclareMathOperator{\SO}{SO}%

\newcommand{\CC}{\mathscr{C}}

\newcommand{\Nscr}{\mathscr{N}}
\newcommand{\Pscr}{\mathscr{P}}

\newcommand{\Dcal}{\mathcal{D}}
\newcommand{\Ocal}{\mathcal{O}}
\newcommand{\Fcal}{\mathcal{F}}

\begin{document}

\title{Counting Minimal Surfaces in Quasi-Fuchsian three-Manifolds}

\author{Zheng Huang}
\address{Department of Mathematics, The City University of New York, Staten Island, NY 10314, USA}
\address{The Graduate Center, The City University of New York, 365 Fifth Ave., New York, NY 10016, USA}
\email{zheng.huang@csi.cuny.edu}

\author{Biao Wang}
\address{Department of Mathematics, Wesleyan University, Middletown, CT 06459, USA}
\email{bwang@wesleyan.edu}


\date{July 17, 2012}


\begin{abstract}
It is well known that every {\qfm} admits at least one closed incompressible {\ms}, and at most
finitely many of them. In this paper, for any prescribed integer $N > 0$, we construct a {\qfm} which
contains at least $2^N$ such {\ms}s. As a consequence, there exists some simple closed {\jc} on
$S_{\infty}^2$ such that there are at least $2^N$ disk-type complete {\ms}s in $\H^3$ sharing this {\jc} as the {\ab}.
\end{abstract}

\maketitle

\section{Introduction}\label{sec:introduction}

Let $M^3$ be a {\htm}, then we can write $M^3 = \H^3/\Gamma$, where $\H^3$ is hyperbolic three-space, and
$\Gamma$ is a discrete subgroup of $\PSL_2(\C)$. The geometry of $M^3$ is largely determined by the behavior
of the action of the discrete group $\Gamma$ on $\mathbb{S}_{\infty}^{2}=\partial \H^3$. When the limit set of
$\Gamma$ is a round circle $\mathbb{S}^1 \subset \mathbb{S}_{\infty}^{2}$, then the corresponding {\tm} $M^3$ is
called {\it Fuchsian}. In this case, $M^3$ contains a (unique) totally geodesic surface $\Sigma$ and $M^3$ is a
warped product space $M^3 = \Sigma \times \R$. The moduli space of Fuchsian manifolds (denoted by
$\Fcal_g(\Sigma)$) is isometric to {\TS}, the space of {\hym}s on the surface $\Sigma$ under the equivalence
induced by the group of diffeomorphisms isotopic to the identity. This space $\Fcal_g(\Sigma)$ is a manifold of
dimension $6g-6$, where $g \ge 2$ is the genus of $\Sigma$. In the case when the limit set of $\Gamma$ lies in a
Jordan curve, $M^3 = \H^3/\Gamma$ is called {\it {\qf}}. This is an important class of {\htm}s and $M^3$ is
diffeomorphic to $\Sigma \times \R$, where $\Sigma$ is a closed surface. The celebrated Simultaneous Uniformization
Theorem of Bers (\cite{Ber72}) implies that each {\qf} manifold is determined by a pair of {\RS}s in {\TS}.  The
space of {\qf} manifolds, called the {\qf} space, is diffeomorphic to the product of {\TS}s, hence of dimension $12g-12$.

{\it All surfaces in this paper, when referred to be contained in a {\qf} manifold, are smooth, closed, incompressible and
of genus greater than one}.

A closed surface embedded in a {\htm} is called {\it incompressible} if it induces a {\fg} injection. Incompressible surfaces
play fundamental roles in the study of {\htm}s (\cite{Has05, Rub07}), and among these surfaces, minimal ones often
provide important geometric information about the ambient {\tm}s (see for instance \cite{Rub05}). Analogs of these surfaces
and {\htm}s are also central objects to study in anti-de-Sitter geometry (see \cite{KS07} and many others).

In this paper, we investigate counting problems on {\ms}s in {\qf} manifolds. We quickly describe our goal of this
work: Any counting problem starts with an (non)existence theorem. It is a fundamental fact, proved by Schoen-Yau
(\cite{SY79}) and Sacks-Uhlenbeck (\cite{SU82}), that any {\qf} manifold $M^3$ contains an area-minimizing immersed
incompressible surface and this surface is shown to be embedded by Freedman-Hass-Scott (\cite{FHS83}). Uhlenbeck
(\cite{Uhl83}) further considered a subclass of {\qf} manifolds which admit a {\ms} of {\pc}s lying in the interval $(-1,1)$,
and such {\qf} manifolds are later (\cite{KS07}) called {\it {\af}}. She showed, among other results, that any {\af} manifold
contains exactly one {\ms}. Therefore one can parametrize the space of {\af} manifolds by studying these {\ms}s (see also
\cite{Tau04}), and obtain important information on {\af} manifolds such as volume of the {\cc} and {\Hd} of the limit set (see \cite{HW11}),
as well as relate these problems to natural metrics on {\TS} (\cite{GHW10}). It was expected (see for instance \cite{Uhl83} page
7) that any (fixed) {\qf} manifold contains at most finitely many (incompressible) {\ms}s. This was indeed the case for stable ones as
shown by Anderson (\cite{And83}) by the method of geometric measure theory. One of the authors here (\cite{Wan12})
showed that there exist {\qfm}s which admit more than one {\ms}. The purpose of this paper is to construct {\qfm}s
that admit arbitrarily many {\ms}s, or more precisely, to prove the following main result:

\bt\label{main}
For any given positive integer $N$, there exists a {\qf} manifold $M^3$ that contains at least $2^N$ distinct (closed and
incompressible) {\ms}s, and each of them is embedded in $M^3$.
\et

We remark here that for different integers $N_1$ and $N_2$, the {\ms}s in the {\qf} manifolds from our construction
might have different genera. In fact, as $N$ gets large, the genera of the {\ms}s are expected to get very large as well.
There is also a different type of counting problem for closed {\ms}s in {\htm}s studied in \cite{HL12}.

When we lift this {\qf} manifold to $\H^3$, a direct consequence of Theorem ~\ref{main} is that we obtain a {\jc}
(not smooth since it is the limit set of a {\qf} group) which bounds many disk-type complete {\ms}s, namely,

\bt
For any given positive integer $N$, there exists a Jordan curve $\Lambda$ on $S_{\infty}^{2}$, such that there exist at least
$2^N$ distinct complete {\ms}s in $\H^3$, sharing the boundary $\Lambda$ at the infinity.
\et

This is related to a type of ``{\it asymptotic Plateau problem}" studied by Morrey (\cite{Mor48}), Almgren-Simon (\cite{AS79})
and Anderson (\cite{And83}) via geometric measure theory, and more recently Coskununzer (\cite{Cos09}) from a
more topological approach. Our result is restricted to $n = 3$, which is very different from higher dimensions. Moreover, many
results on this asymptotic Plateau problem concern the regularity of the prescribed simple closed curve at infinity (see for
example \cite{HL87,GS00} and many others), the {\jc} in our construction is the limit set of some {\qfg}, therefore, except in
the Fuchsian case, this curve is well-known to contain no rectifiable arcs (see for instance \cite{Leh64, Ber72}).


\subsection*{Outline of the construction}
Our construction can be outlined as follows: we separate the unit ball $\B^3$ into $N+1$ chambers by $N$ parallel Euclidean
circles and these circles are slightly fattened to form bands (see figure 3); bands are then connected by narrow bridges
to form one piecewise smooth {\jc} $\Lambda$ (see figure 4) on the sphere at infinity $S_{\infty}^{2}$. There will be
{\tdc} next to every bridge to form a {\mct}; we then cover $\Gamma$ by an even number of small circles on $S_{\infty}^{2}$, and
they induce, via inversions (or reflections), a {\qfg} $\Gamma_\Lambda$ whose limit set is within a small neighborhood of the {\jc} $\Lambda$,
and hence this gives rise to a {\qfm} (notably this construction of {\qfg}s was known to Poincar\'e and Fricke-Klein); at last,
there are more than $2^N$ distinct ways to arrange pairs of circles next to the bridges to bound {\mct}s (see figure 6), and for
each arrangement, we obtain a compact region of mean convex boundary inside the {\cc} of the corresponding {\qfm}, which
allows us to apply the results of Schoen-Yau and Sacks-Unlenbeck to trap a {\ms} in this region.

There is a construction (unpublished) due to Hass-Thurston of {\qfm}s which admit many {\ms}s (see \cite{GW07}).
Their construction works for all genera greater than one, while our construction is quite different, and of quantitative nature.
Moreover, one can see in our construction that each {\ms} is obtained by removing some solid {\mct}s in the regions within
the {\cc} of the {\qfm}. We (\cite{HW11}) have used a simplified version of this type of construction to obtain a {\qf} manifold
which does not admit a foliation of closed surfaces of {\cmc}.


\subsection*{Plan of the paper}
We organize this paper as follows: Subsections \S 2.1, \S 2.2 and \S 2.3 consist of preliminaries where we recall {\qfg}s, {\qfm}s,
{\ms}s in {\qfm}s and {\ms}s of catenoid type in $\H^3$. We prove the main Theorem ~\ref{main} in Section \S 3 via a construction
outlined above. There are several subsections in this main section: in \S 3.1, we obtain a condition when {\tdc} on $S_{\infty}^2$
bound a {\mct}; in \S 3.2, we follow the classical construction of {\qfg}s by using inversions of circles which bound open disks
covering some piecewise smooth {\jc}; in \S 3.3, we find a {\ms} in the corresponding {\qfm} in the region of the {\cc} with some
solid catenoid removed; the subsection \S 3.4 contains a topological lemma which will be used later to show the {\ms}s obtained 
from the construction are distinct; in \S 3.5, we assemble these results to finalize the proof.


\subsection*{Acknowledgements} The idea of using {\mct}s as barrier surfaces was originally suggested by
Bill Thurston. We dedicate this paper to his memory. We are also grateful to Joel Hass for explaining the 
Hass-Thurston construction on {\qfm}s that contain arbitrarily many {\ms}s. We would also like to express 
our sincere gratitude to the referee for numerous corrections, many constructive suggestions and helpful 
comments. The research of Z. H. is partially supported by a CUNY-CIRG award $\#1861$, and Provost's 
Research Scholarship of CUNY-CSI.


\section{Preliminaries}
In this section, we collect some important facts and describe some key properties on {\qf} groups, {\qf} manifolds
and their {\ms}s.
\subsection{Quasi-Fuchsian groups and {\qfm}s}

In this paper, we will work in the ball model ($\B^3$) of the hyperbolic three-space $\H^{3}$, i.e.,
\begin{equation*}
   \B^{3}=\{(x,y,z)\in\R^{3}\ |\
   x^{2}+y^{2}+z^{2}<{}1\},
\end{equation*}
equipped with metric
\begin{equation*}
   ds^{2}=\frac{4(dx^{2}+dy^{2}+dz^{2})}{(1-r^{2})^{2}}\ ,
\end{equation*}
where $r=\sqrt{x^{2}+y^{2}+z^{2}}$. The hyperbolic space $\H^{3}$ has a natural compactification:
$\overline{\H^{3}}=\H^{3}\cup{}S_{\infty}^{2}$, where $S_{\infty}^{2}\cong\C\cup\{\infty\}$ is the Riemann
sphere. The orientation preserving isometry group of the three-ball $\B^3$ is denoted by $\Mob(\B^3)$, which
consists of M\"obius transformations that preserve the unit ball (see \cite[Theorem 1.7]{MT98}). It's well known
that $\Mob(\B^3)\cong\Isom^+(\H^3)$, which is isomorphic to $\PSL_{2}(\C)$.

Suppose that $X$ is a subset of $\B^{3}$, we define the {\it {\ab}} of $X$ by
\beq
   \partial_{\infty}X=\overline{X} \cap{}S_{\infty}^{2}\ ,
\eeq
where $\overline{X}$ is the closure of $X$ in $\overline{\H}{}^{3}$.

Using the above notation, if $P$ is a geodesic plane in $\B^{3}$, then $P$ is perpendicular to $S_{\infty}^{2}$
and its {\ab} $C\stackrel{\text{def}}{=}\partial_{\infty}P$ is an Euclidean circle on $S_{\infty}^{2}$. We also say that
$P$ is {\it asymptotic to} $C$.

A (torsion free) discrete subgroup $\Gamma$ of $\Isom^+(\H^{3})$ is called a {\it {\kg}}, and
\beq
   M^3 \equiv{}M_{\Gamma}^3=\H^{3}/\Gamma
\eeq
is a complete {\htm} with the {\fg} $\pi_{1}(M^3)\cong\Gamma$.

For any Kleinian group $\Gamma$, $\forall\,p\in\H^{3}$, its orbit set
\beq
   \Gamma(p)=\{g(p)\ |\ g\in{}\Gamma\}
\eeq
has accumulation points on the Riemann sphere $S_{\infty}^{2}$, which are called the {\it limit points} of $\Gamma$.
The {\it limit set} of $\Gamma$, denoted by $\Lambda_{\Gamma}$, is the closure of the limit points on
$S_{\infty}^{2}$, i.e. $\Lambda_\Gamma=\overline{\Gamma(p)}\cap S_{\infty}^{2}$. It is known that
$\Lambda_\Gamma$ is independent of the choice
of the reference point $p$, and it is a closed $\Gamma$-invariant subset of $S_{\infty}^{2}$. The open set
$\Omega_{\Gamma}=S_{\infty}^{2}\setminus\Lambda_{\Gamma}$ is called the {\it domain of discontinuity}.
The {\kg} $\Gamma$ acts properly discontinuously on $\Omega_\Gamma$, and the quotient
$\Omega_\Gamma/\Gamma$ is a finite union of {\RS}s of finite type if $\Gamma$ is finitely generated
(see \cite{Mar74}).

For any {\kg} $\Gamma$, the {\it convex hull} of the limit set $\Lambda_\Gamma$,
denoted by $\Hull(\Lambda_{\Gamma})$, is the smallest convex subset in $\H^{3}$ whose closure in
$\H^{3}\cup S_{\infty}^{2}$ contains $\Lambda_\Gamma$. The quotient space $
\CC_{\Gamma}=\Hull(\Lambda_\Gamma)/\Gamma$ is called the {\it {\cc}} of $M^3$.

A {\kg} $\Gamma$ is {\it {\qf}} if $S_{\infty}^{2}\setminus\Lambda_{\Gamma}$ has exactly two components, denoted
by $\Omega_{\pm}$, such that each component is invariant under $\Gamma$. If $\Gamma$ is a {\qf} group, then
its limit set $\Lambda_{\Gamma}$ is a Jordan curve, and the quotient space $M=\H^3/\Gamma$, which is
called a {\it {\qfm}}, is diffeomorphic to $\Sigma \times \R$, where $\Sigma=\Omega_{+}/\Gamma$ or
$\Omega_{-}/\Gamma$ is a finitely punctured compact surface (see \cite[Lemma 3.2 and 3.3]{Mar74}).

We always assume that the {\qf} group $\Gamma$ is {\it torsion-free}. Therefore $\Omega_{+}/\Gamma$ and
$\Omega_{-}/\Gamma$ are closed surfaces.

\subsection{Minimal surfaces in {\qfm}s}
Suppose that $\Sigma$ is a surface (compact or complete) and that $M$ is a $3$-dimensional Riemannian
manifold (compact or complete). An immersion $f:\Sigma\to{}M$ is called a {\it {\ms}} if its {\mc} is identically
equal to zero.

If $\Sigma$ is a closed {\ms}, then it is of {\it least area} if its area is less than that of any other surfaces in the
same homotopy class; it is {\it area minimizing} if its area is no larger than that of any surface in the same
homology class. If $\partial\Sigma\ne\emptyset$, we require that other surfaces should share the same
boundary as $\Sigma$. If $\Sigma$ is a (non-compact) complete minimal surface, then it is
{\it least area} or {\it area minimizing} if any compact subdomain $K$ of $\Sigma$ is least area or area minimizing.

Recall that a map $f: \Sigma \to{}M$ of a closed surface $\Sigma$ into a $3$-manifold $M$ is
called {\it incompressible} if $f_*:\pi_1(\Sigma)\to\pi_1(M)$ is injective.

The {\cc} of the {\qf} group $\Gamma$ is a compact {\htm}, whose boundary is convex {\wrt} the inward normal
vector. Schoen-Yau (\cite{SY79}) and Sacks-Uhlenbeck (\cite{SU82}) proved that such a {\htm} always contains
an (immersed) incompressible least area {\ms}. By the work of Freedman, Hass and Scott
(\cite[Theorem 5.1]{FHS83}), this {\ms} is actually embedded. Here we should remark that the ambient {\htm}
appeared in the main results in \cite{SY79,SU82,FHS83} are stated to be compact without boundary, i.e., closed.
But from the discussions in \cite[p. 635]{FHS83} and the works of Meeks-Yau (\cite{MY82a, MY82b}),
these results can be extended to the compact {\htm} whose boundaries are {\it convex} or {\it mean convex}.


\subsection{Minimal catenoids in $\H^3$} Minimal surfaces of catenoid type in $\H^3$ whose {\ab} consists of two
circles on $S_{\infty}^{2}$ will play an important role in the construction: they serve as barriers among {\ms}s. In this
subsection, we describe the existence of these minimal catenoids. The properties of minimal catenoids can be found
in many articles, for example, \cite{Mor81,Gom87,BSE10,Seo11,Wanar}.

Let $G$ be the subgroup of $\Mob(\B^3)$ that leaves a geodesic $\gamma\subset\B^{3}$ pointwisely fixed, then we
call $G$ the {\it spherical group} of $\B^{3}$ and $\gamma$ the {\it rotation axis} of $G$. A (connected) surface $\Pi$
in $\B^{3}$ that is invariant under $G$ is called a {\it spherical surface} or a {\it surface of revolution}. If $\Pi$ is
minimal, then it is called a {\it spherical minimal catenoid} or a {\it minimal catenoid}.

For two circles $C_{1}$ and $C_{2}$ in $\H^{3}$, if there is a geodesic $\gamma$ such that each of the circles
$C_{1}$ and $C_{2}$ is invariant under the group of rotations that fixes $\gamma$ pointwisely, then $C_{1}$ and
$C_{2}$ are said to be {\it coaxial}, and $\gamma$ is called the {\it rotation axis} of $C_{1}$ and $C_{2}$. If $C_{1}$
and $C_{2}$ are two disjoint circles on $S_{\infty}^{2}$, then they are always coaxial. In this case, we want to know
whether there exists a spherical minimal surface that is asymptotic to $C_{1}\cup{}C_{2}$. This existence depends
on two kinds of distances between two disjoint circles on $S_{\infty}^{2}$, which we now describe.

Equipping the Riemann sphere $S_\infty^2$ with the spherical metric $\rho$, then every geodesic on $(S_\infty^2,\rho)$ 
is just a great circle. The orientation preserving isometry group of $S_\infty^2$, denoted by $\Isom^+(S_\infty^2)$, is 
$\SO(3)$ which is a subgroup of $\Mob(\B^3)$. Suppose $C_{1}$ and $C_{2}$ are disjoint (round) circles in $S_{\infty}^{2}$ 
which bound disjoint subdisks $\Delta_1$ and $\Delta_2$ of $S_{\infty}^{2}$ with injective radii $\leq\pi/2$. For $i = 1,2$, 
let $P_{i}$ be the geodesic plane in $\B^3$ such that $\partial_{\infty}P_{i}=C_{i}$. Then we may define two distances 
between $C_{1}$ and $C_{2}$ by
\begin{align}
   \label{d_L}
   d_{L}(C_{1},C_{2})&=\dist(P_{1},P_{2})\ ,\\
   \label{rho}
   \rho(C_{1},C_{2})&=\min\{\rho(p,q)\ |\ p\in{}C_{1}\
      \text{and}\ q\in{}C_{2}\}
\end{align}
where $\dist(\cdot,\cdot)$ is the hyperbolic distance in $\B^3$ (or $\H^3$).
The distance $d_{L}$ is invariant under $\Mob(\B^3)$, whereas the spherical
distance $\rho$ is only invariant under $\SO(3)$.

Suppose that $G$ is the spherical group of $\B^{3}$ with rotation axis $\gamma_{0}$, a geodesic along on the $y$-axis,
then $\B^3/G\cong\B_{+}^2$, where
\begin{equation}
   \B_{+}^2=\{(x,y,z)\in\B^3\ |\ z = 0, \ y\geq{}0\}\ .
\end{equation}
If $\Pi$ is a spherical {\mct} in $\B^3$ {\wrt} the axis $\gamma_{0}$, then the curve $\sigma=\Pi\cap\B_{+}^{2}$ is
called the {\it generating curve} of $\Pi$. Gomes (\cite{Gom87}) proved that $\sigma$ is symmetric about the
$y$-axis up to isometries. Moreover, he proved the following theorem:

\bt(\cite[Proposition 3.2]{Gom87})\label{Gomes}
There exists a finite constant $d_{0}>0$ such that for two disjoint circles $C_{1},C_{2}\subset{}S_{\infty}^{2}$, if
$d_{L}(C_{1},C_{2})\leq{}d_{0}$, then there exist a {\ms} $\Pi$ which is a surface of revolution asymptotic to
$C_{1}\cup{}C_{2}$.
\et

In this paper, we require the existence of {\it least area} minimal catenoids in order to apply the results in
\cite{MY82b}. This is fulfilled by:

\bt(\cite[Theorem 1.2]{Wanar})\label{wang}
There exists a finite constant $d_{1}>0$ such that for {\tdc} $C_{1},C_{2}\subset{}S_{\infty}^{2}$, if
$d_{L}(C_{1},C_{2})\leq{}d_{1}$, then there exist a least area {\ms} $\Pi$ which is a surface of revolution
asymptotic to $C_{1}\cup{}C_{2}$.
\et


\section{Main Construction}

This section is devoted to constructing a {\qfm} which admits at least $2^N$ {\ms}s, for any prescribed positive
integer $N$. Let us outline the organization of this section. In \S 3.1, in order to utilize the Theorems
~\ref{Gomes}, ~\ref{wang}, to find a {\mct} asymptotic to two disjoint circles on $S_\infty^2$, we show know how to compute
the distance $d_L$ (see the equation \eqref{d_L}) between two circles (Lemma ~\ref{0dist}); in \S 3.2, we use
inversions {\wrt} disjoint circles on $S_\infty^2$ to construct a {\qfg} whose limit set is contained within a small
neighborhood of a prescribed (piecewise smooth) Jordan curve on $S_\infty^2$ (Theorem ~\ref{inversion}); in
\S 3.3, we use the notion of {\fp} (Definition ~\ref{fp}) for a {\qfg} to find a compact {\htm} with mean convex
boundary, which allows us to apply fundamental results of Schoen-Yau, Saks-Uhlenbeck to find a closed {\ms}
(Theorem ~\ref{key}); finally in \S 3.4, we combine the results in previous subsections to prove our main theorem.

\subsection{Distances between circles on $S_\infty^2$} We have introduced two distances (see equations
\eqref{d_L} and \eqref{rho}) between two disjoint circles on $S_\infty^2$. The purpose of this subsection is to find
the relationship between these two distances, and use it to find conditionsunder which two circles determine a {\mct}.

\begin{figure}[htbp]
\begin{minipage}[t]{0.75\linewidth}
\begin{center}
   \includegraphics[scale=1]{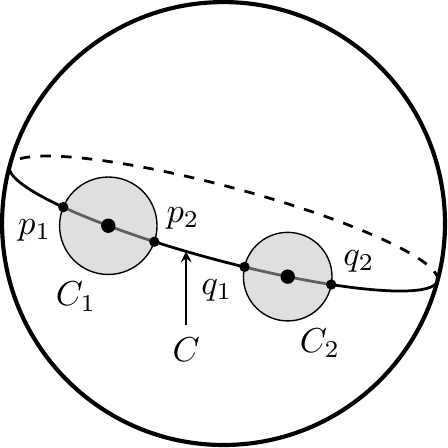}
\end{center}
\end{minipage}
\caption{Spherical distance between two circles}\label{fig:spherical distance}
\end{figure}

We may calculate the above two distances as follows. As in the subsection \S 2.3, let $C_{1}$ and $C_{2}$ be {\tdc} on
$S_{\infty}^{2}$, and $P_{i}$ be the geodesic planes in $\B^3$ and $\Delta_{i}$ be the disks on $S_\infty^2$ with
injective radius $\leq\pi/2$ such that $\partial_{\infty}P_{i}=C_{i}$ and $\partial{}\Delta_{i}=C_{i}$ for $i=1,2$ and such
that $\bar{\Delta}_1 \cap \bar{\Delta}_2 = \emptyset$. Now let
$C$ be the great circle that passes the centers of $\Delta_1$ and $\Delta_2$, and we also mark $p_1$ and $p_2$ as
the intersection points of $C_1$ and $C$, while $q_1$ and $q_2$ as the intersection points of $C_2$ and $C$. See
Figure \ref{fig:spherical distance}.

Let $P$ be the geodesic plane asymptotic to $C$ and $\gamma_i=P_i\cap{}P$. Then from the definitions, we have
$d_{L}(C_{1},C_{2})=\dist(\gamma_{1},\gamma_{2})$, and the spherical distance $\rho(C_{1},C_{2})$ is equal to the
shorter length of two geodesic segments on $C$ between $\Delta_1$ and $\Delta_2$. One of our key lemmas is the
following:

\bl\label{0dist}
Suppose that $C_{1}$ and $C_{2}$ are t{\tdc} on $S_{\infty}^{2}$ such that
$\overline\Delta_{1}\cap{}\overline\Delta_{2}=\emptyset$, where $\Delta_i\subset{}S_\infty^2$ is a (spherical) disk
bounded by $C_i$ with injective radius $\leq\pi/2$. Then
\begin{equation*}
   d_{L}(C_{1},C_{2})\to{}0\ ,
   \quad\text{as}\ \rho(C_{1},C_{2})\to{}0\ .
\end{equation*}
In particular, if $C_{1}$ and $C_{2}$ are tangent at one point on
$S_{\infty}^{2}$, then $d_{L}(C_{1},C_{2})=0$.
\el

We note, in this lemma, $\rho(C_{1},C_{2})\to{}0$ means that the spherical distance between two closed sets $C_1$
and $C_2$ approaches zero while the injective radii of two disks $\Delta_1$ and $\Delta_2$ stay unchanged.

\bp
Without lost of generality,
we may assume that the center of the disk $\Delta_1$ is the north pole, i.e. $(0,0,1)$, and the center of
$\Delta_2$ lies on the great circle $C=\{(x,0,z)\in\R^3 \ |x^2+z^2=1\}$. In fact, using at most two rotations along
some axes, we can obtain the required picture. Since these two rotations are elements of $\SO(3)$, the two distances
\eqref{d_L} and \eqref{rho} are preserved.

\begin{figure}[htbp]
\begin{minipage}[t]{0.75\linewidth}
\begin{center} 
   \includegraphics[scale=1]{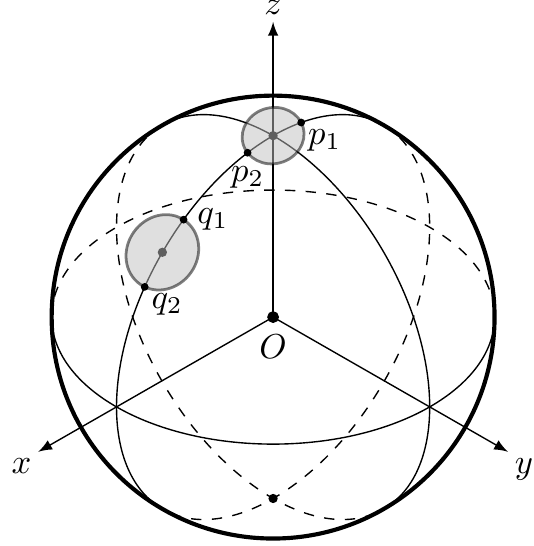}
\end{center}
\end{minipage}
\caption{Stereographic projection of two circles}\label{fig:two distances}
\end{figure}

As above, we suppose that $C_1\cap{}C=\{p_1,p_2\}$ and $C_2\cap{}C=\{q_1,q_2\}$. We order them on the great
circle $C$ as $p_1$, $p_2$, $q_1$ and $q_2$, see Figure \ref{fig:two distances}.

Then it is easy to see that $\rho(C_{1},C_{2})=\min\{\rho(p_2,q_1),\rho(p_1,q_2)\}$. We may assume that
$\rho(C_{1},C_{2})=\rho(p_2,q_1)$.

For $i=1,2$, let $\gamma_i=P_i\cap{}P$, where $P=\{(x,0,z)\in\R^3 \ |x^2+z^2\leq{}1\}$, and $P_i$ is the geodesic
plane asymptotic to $C_i$. Then $\partial_\infty\gamma_1=\{p_1,p_2\}$ and
$\partial_\infty\gamma_2=\{q_1,q_2\}$. Since $\dist(P_{1},P_{2})=\dist(\gamma_1,\gamma_2)$, we have
$d_{L}(C_{1},C_{2})=\dist(\gamma_1,\gamma_2)$.

There exists a M\"obius transformation $\Phi\in\Mob(\widehat{\R}^3)$ that satisfies the following properties
(see \cite[pp.19-20]{MT98}):
\begin{itemize}
   \item $\Phi|{\B^3}:\B^3\to\H^3$ is an isometry, and
   \item $\Phi|{S_\infty^2\setminus\{e_3\}}:S_\infty^2\setminus\{e_3\}
         \to\C$ is the stereographic projection,
\end{itemize}
where $\widehat{\R}^3=\R^3\cup\{\infty\}$ and $e_3=(0,0,1)$.

For $i=1,2$, let $\tilde\gamma_i=\Phi(\gamma_i)$, and let $u_i=\Phi(p_i)$ and $v_i=\Phi(q_i)$, respectively. Then
$\tilde\gamma_1$ and $\tilde\gamma_2$ are two geodesics in the $xy$-plane. Since $\Phi$ preserves the
{\hym}s, we have
\be\label{eq-dist}
   \dist(\tilde\gamma_1,\tilde\gamma_2)=
   \dist(\gamma_1,\gamma_2)=d_{L}(C_1,C_2)\ .
\ene
Suppose the points $u_1$, $u_2$, $v_1$ and $v_2$ are contained on the $x$-axis in the following order:
$u_1,v_2,v_1,u_2$. Since $\partial_\infty\tilde\gamma_1=\{u_1,u_2\}$ and
$\partial_\infty\tilde\gamma_2=\{v_1,v_2\}$, applying \cite[Eq. (7.23.1)]{Bea95}, we have
\beq
   \tanh^2[\dist(\tilde\gamma_1,\tilde\gamma_2)]=
   \frac{1}{[u_1,v_2,v_1,u_2]}=
   \frac{(u_1-v_2)(v_1-u_2)}{(u_1-v_1)(v_2-u_2)}\ .
\eeq
Then we have the following consequences:
\begin{equation*}
    \rho(C_1,C_2)\to{}0\
    \ \Longrightarrow\
    q_1\to{}p_2
    \ \Longrightarrow\
    v_1\to{}u_2 \\
    \ \Longrightarrow\
    \dist(\tilde\gamma_1,\tilde\gamma_2)\to{}0\ ,
\end{equation*}
and then $d_{L}(C_{1},C_{2})\to{}0$ because of \eqref{eq-dist}.
\ep
This lemma enables us to conclude that if we can arrange two disjoint circles on $S_{\infty}^{2}$ sufficiently
close {\wrt} the (spherical) distance $\rho$, then we can apply the Theorem ~\ref{Gomes} (and the
Theorem ~\ref{wang}) to find a (least area) {\mct} asymptotic to these two circles.

\subsection{Constructing {\qfg}s via inversions} In this subsection, we use inversions of circles on
$S_\infty^2$ to generate a {\qfg} and we obtain a {\qfm} in the quotient.

Suppose that $C_1$ and $C_2$ are two circles on $S_\infty^2$. Let $f_i$ be the inversions {\wrt} to the
circle $C_i$, then each inversion is of order $2$, i.e. $f_{i}^{2}=\id$, for $i=1,2$. Consider the M\"obious
transformation $g=f_1\circ{}f_2$. We have:
\bl\label{commute}
For this element $g$, there are three possibilities:
\begin{enumerate}
   \item
   If $C_1\cap{}C_2=\emptyset$, then $g$ is a loxodromic element, whose fixed points are contained in
   the disks $\Delta_1$ and $\Delta_2$, respectively.
   \item If $C_1$ and $C_2$ are tangent at some point $z_0\in\R^2$, then $g$ is a
         parabolic element with exactly one fixed point $z_{0}$.
   \item If $C_1$ and $C_2$ are orthogonal, and $C_1\cap{}C_2=\{z_1,z_2\}$, then $g$ is an elliptic element
         of order $2$ with fixed points $\{z_1,z_2\}$.
\end{enumerate}
\el
\bp
We will only need to prove the third case, for which we work in the upper half space model. Suppose $C_1$
is a circle with radius $r_1$ whose center is the origin, and $C_2$ is a circle of radius $r_2$ whose center
is at the point $(a,0)$, where
\beq
   a=(r_{1}^{2}+r_{2}^{2})^{1/2}\ .
\eeq
For any $x\in\R^2$, let $x^*=x/|x|^2$, where $|x|$ is the Euclidean norm. Then the inversions with respect to
$C_1$ and $C_2$ have the following forms:
\beq
   f_1(x)=r_{1}^{2}x^{*}
   \quad\text{and}\quad
   f_2(x)=a+r_{2}^{2}(x-a)^{*}\ .
\eeq
We need show that these two inversions commute, namely, $f_1\circ{}f_2=f_2\circ{}f_1$. It is clear that
\be\label{f2f1}
   f_2\circ{}f_1(x)=a+r_{2}^{2}(r_{1}^{2}x^{*}-a)^{*}\ .
\ene
Direct computation shows
\begin{align}
   |f_2(x)|^2&=\frac{|r_{1}^{2}x^{*}-a|^{2}|x|^{2}}{|x-a|^{2}},
   \end{align}
   and
   \begin{align}
   |r_{1}^{2}x^{*}-a|^2-r_{2}^{2}&=\frac{r_{1}^{2}(|x-a|^2-r_2^2)}{|x|^2}\ .
\end{align}
Therefore, we have
\begin{align*}
    f_1\circ{}f_2(x)
       &=\frac{r_{1}^{2}|x-a|^{2}}{|r_{1}^{2}x^{*}-a|^{2}|x|^{2}}
         \left(a+r_{2}^{2}(x-a)^{*}\right)\\
       &=\frac{1}{|r_{1}^{2}x^{*}-a|^{2}}
         \left(\frac{r_{1}^{2}|x-a|^{2}a+r_{1}^{2}r_{2}^{2}(x-a)}{|x|^2}\right)\\
       &=\frac{1}{|r_{1}^{2}x^{*}-a|^{2}}
         \left(\frac{r_{1}^{2}(|x-a|^2-r_{2}^{2})}{|x|^{2}}\,a+
         r_{1}^{2}r_{2}^{2}x^{*}\right)\\
       &=\frac{1}{|r_{1}^{2}x^{*}-a|^{2}}
         \left((|r_{1}^{2}x^{*}-a|^2-r_{2}^{2})a+r_{1}^{2}r_{2}^{2}x^{*}\right)\\
       &=a+r_{2}^{2}(r_{1}^{2}x^{*}-a)^{*}\ .
\end{align*}
This agrees with \eqref{f2f1}. Since that $f_{i}^2=\id$ for $i=1,2$, we have
\beq
g^2=f_1\circ{}f_{2}\circ{}f_{1}\circ{}f_{2}=f_1\circ{}f_{1}\circ{}f_{2}\circ{}f_{2}=\id.
\eeq
Recall that any inversion {\wrt} some circle fixes that circle pointwisely, then it is not hard to see that
$\{z_1,z_2\}$ are the points that are fixed by both $f_1$ and $f_2$. This completes the proof of the lemma.
\ep
For any piecewisely smooth simple closed curve $\Lambda$ on $S_\infty^2$, we can cover it by finitely many small disks.
the above lemma allows us to construct a {\qfg} $\Gamma$, whose limit set $\Lambda_\Gamma$
is within some small distance of the curve $\Lambda$.

We now proceed with this construction as follow: let
$\Nscr_\delta(\Lambda)\subset{}S_\infty^2$ be a $\delta$-neighborhood of $\Lambda$, where $\delta$ is a
small positive number. Then we can cover the curve $\Lambda$ by a family of finitely many open disks
$\{\Delta_i\}_{i=1}^{2L}$ on $S_\infty^2$, where $L\geq{}3$, which satisfy the following conditions:
\begin{itemize}
   \item For each disk $\Delta_i$, its center is on $\Lambda$
         and its injective radius is $\leq\delta$.
   \item Let $C_i=\partial{}\Delta_i$, each circle
         $C_i$ intersects $C_{i+1}$ at an angle of $\pi/2$ and
         $C_i\cap{}C_j=\emptyset$ for $|i-j|\geq{}2$, here
         $C_{2L+1}\stackrel{\text{def}}{=}C_1$.
   \item Suppose that $C_{i}\cap{}C_{i+1}=\{z_i,z_i'\}$, $i=1,\ldots,2L$.
\end{itemize}
It is easy to see that we have
$\Lambda\subset\bigcup\limits_{i=1}^{2L}\Delta_i\subset
\Nscr_\delta(\Lambda)$.

Let $f_{i}$ be the inversion with respect to the circle $C_{i}$, and let $g_{i}=f_{i}\circ{}f_{i+1}$, $i=1,2,\ldots,2L$,
where $C_{2L+1}=C_1$. Let $\Gamma'$ be the group generated by the inversions $f_1,f_2,\ldots,f_{2L}$.
Let $\Gamma_{0}$ be the group whose elements are the product of {\it even number} of inversions $f_{i}$.
Clearly $\Gamma_{0}$ is a subgroup of $\Gamma'$ of index 2. The main result in this subsection is the
following theorem:

\bt\label{inversion}
Let $\Gamma'$ and $\Gamma_0$ be the above groups, and let $\Lambda_{\Gamma'}$ be the limit set of
$\Gamma'$, then we have
\ben
\item
$\Lambda_{\Gamma'}$ is a {\jc} that passes through all of the points $\{z_i,z_i'\}_{i=1}^{2L}$. Moreover
$\Gamma_0$ is a {\qfg} and $\Lambda_{\Gamma_0}=\Lambda_{\Gamma'}\subset \Nscr_\delta(\Lambda)$.
\item
The group $\Gamma_{0}$ has the following presentation
\be
   \Gamma_{0}=\left\langle{}g_{1},\ldots,g_{2L}\ |\
              g_{1}^{2}=1,\,\ldots,\,g_{2L}^2=1,\,
              g_{1}\cdots{}g_{2L}=1\right\rangle\ .
\ene
\item
The quotient space $\Ocal=\Omega_{+}/\Gamma_{0}$ is a genus zero orbifold that has $2L$ cone points with
cone angle $\pi$, where $\Omega_{\pm}=S_\infty^2\setminus\Lambda_{\Gamma_0}$.
\een
\et
\br
This construction of {\qfg}s was known at the end of the 19th century by Poincar\'e and Fricke-Klein
(see \cite{Poi83,Poi85} and \cite[pp. 399-445]{FK65}, see also \cite[p. 263]{Ber72}).
\er
\bp
\noindent
(1). Similar to the discussion in \cite[Chapter IV]{Mag74} or \cite[Chapter 6]{MSW02}, one can show that
$\Lambda_{\Gamma'}$ is a Jordan curve contained in $\cup{}\Delta_{i}\subset\Nscr_\delta(\Lambda)$.
By the construction, $\Gamma_{0}$ is a subgroup of $\Gamma'$ with index 2, so $\Gamma_{0}$ is a normal
subgroup of $\Gamma'$, then it is well known (see for instance \cite[Lemma 2.22]{MT98}) that they have the
same limit set, namely, $\Lambda_{\Gamma_{0}}=\Lambda_{\Gamma'}$.

(2). By the Lemma ~\ref{commute}, we have $g_{1}^{2}=1$, $\ldots$,
$g_{2L}^2=1$, and $g_{1}\cdots{}g_{2L}=f_{1}^2=1$. In order to show that $\Gamma_0$ is generated by
$g_1,\ldots,g_{2L}$, we only need to consider the special case, i.e., suppose that $h=f_{i}f_{j}$, here $|i-j|=k\geq{}2$.
Then we may write $h$ in the form
\begin{equation*}
   h=\begin{cases}
        g_{i}g_{i+1}\cdots{}g_{j-1} & \text{if}\ i<j\ ,\\
        g_{i-1}\cdots{}g_{i-k+1}g_{j} & \text{if}\ i>j\ ,
     \end{cases}
\end{equation*}
which implies that any element of $\Gamma_{0}$ can be written by the
product of finitely many $g_{i}$, $i=1,\ldots,2L$.

(3). The domain $S_{\infty}^{2}\setminus\bigcup\limits_{i=1}^{2L}\Delta_{i}$ has
two components, denoted by $E$ and $E'$, where $E$ is a $2L$-polygon
with vertices $z_1,z_2,\ldots,z_{2L}$ and $E'$ is also a
$2L$-polygon with vertices $z_{1}',z_{2}',\ldots,z_{2L}'$. Both $E$
and $E'$ are orbifolds that contain $2L$ corner points with corner
angle $\pi/2$. Each of them is a fundamental domain of $\Gamma'$ in
$S_\infty^2$. Recall that $\Gamma_{0}$ is a subgroup of $\Gamma'$ with index $2$,
let $\Ocal$ be the double cover of $E$ or $E'$, then $\Gamma_{0}$ is the
fundamental group of $\Ocal$ (see \cite[pp. 423-424]{Sco83}) and
$\Omega_{+}/\Gamma_{0}=\Ocal$ is a $2$-sphere that has $2L$ cone points
with cone angle $\pi$.
\ep

It is well known that $\Ocal$ can be doubly covered by a closed surface $\Sigma$ of genus
$g=L-1 \geq 2$. Therefore $\Gamma_{0}$ contains a torsion-free subgroup $\Gamma$ of
index $2$. Therefore again this subgroup is normal and we have
$\Lambda_{\Gamma}=\Lambda_{\Gamma_0}$, and $\H^3/\Gamma\cong{}\Sigma\times\R$ is
{\qf}. In particular, we have showed:

\begin{cor}
Suppose that $\Lambda$ is a simple closed (piecewise smooth) curve on $S_\infty^2$, then for any positive number
$\delta$, there exists a torsion-free {\qfg} $\Gamma$ such that its limit set $\Lambda_\Gamma$ is contained in a
$\delta$-neighborhood of $\Lambda$. Furthermore, suppose that $\Sigma$ is a finite type surface such that
$\H^3/\Gamma\cong{}\Sigma\times\R$, then $\genus(\Sigma)=g(\Lambda,\delta)$, here $g(\Lambda,\delta)$ is a positive integer
depending only on $\Lambda$ and $\delta$.
\end{cor}


\subsection{Fundamental polyhedron}
The goal of this subsection is to obtain a compact {\htm} $Y$ whose boundary is mean convex
{\wrt} the inward normal vector. This is achieved by arranging disjoint circles on $S_{\infty}^{2}$ in a
particular way such that we can use {\mct}s (asymptotic to {\tdc}) as barrier surfaces. Let us start
with two definitions.
\begin{definition}
Suppose that $\Lambda_{1},\ldots,\Lambda_{n}$ are a family of disjoint Jordan curves on $S_{\infty}^{2}$,
where $n\geq{}3$, we say that $\Lambda_{1},\ldots,\Lambda_{n}$ are in {\it good position} if each
Jordan curve bounds a topological disk on $S_{\infty}^{2}$ that does not contain any other Jordan curves
in the family.
\end{definition}
First we need a lemma on separation:
\begin{lem}\label{disjoint}
Let $C_{1}$, $C_{2}$ and $C$ be three disjoint circles on $S_{\infty}^{2}$ that are in good position.
Suppose that $\Pi$ is the {\mct} with $\partial_{\infty}\Pi=C_{1}\cup{}C_{2}$, and $P$ is the geodesic
plane with $\partial_{\infty}P=C$, then $\Pi\cap{}P=\emptyset$.
\end{lem}

\bp We may assume that $P$ is on the $xy$-plane, i.e. $P=\{(x,y,0)\ |\ x^2+y^2<1\}$, otherwise we
may find a M\"obius transformation $g\in\Mob(\B^3)$ such that $g(P)$ is on the $xy$-plane (see
\cite[Proposition 1.3]{MT98}). In addition, we assume that $C_{1}$ and $C_{2}$ are above the
$xy$-plane. We now need prove that $\Pi$ is above $P$.

If $\Pi\cap{}P=\emptyset$, then we are done. Otherwise, for each $t\in(-1,0]$, let $P(t)$ be the
geodesic plane that is perpendicular to the $z$-axis at the point $(0,0,t)$. Obviously $P=P(0)$.
It is easy to see that the family of the geodesic planes $\{P(t)\}_{-1<t\leq{}0}$ foliates the
lower half ball, i.e $\{(x,y,z)\in\B^3\ |\ -1<z\leq{}0\}$. If $\Pi\cap{}P\ne\emptyset$, then we always can
find a number $t_{0}\in(-1,0)$ such that $P_{t_{0}}$ is tangent to $\Pi$ at some point in $\B^3$.
But this is impossible by the maximum principle, therefore $P$ and $\Pi$ must be disjoint.
\ep

Suppose that $\Gamma$ is a {\qfg} such that $\B^3/\Gamma\cong{}\Sigma\times\R$, where $\Sigma$
is a closed surface of $\genus\geq{}2$. The limit set of $\Gamma$ is denoted by $\Lambda_\Gamma$,
and $\Lambda_\Gamma$ separates $S_{\infty}^{2}$ into two components, namely,
$S_{\infty}^{2}\setminus\Lambda_\Gamma=\Omega_{+}\cup\Omega_{-}$.
\begin{definition}\label{fp} (see section 2.1.2 in \cite{MT98})
A  closed convex set $\Pscr$ in the hyperbolic space $\B^3$ bounded by finite collection of hyperbolic planes is
called a {\it fundamental polyhedron} for a {\qfg} $\Gamma$, if it satisfies the following conditions:
\begin{enumerate}
   \item $\displaystyle\bigcup_{g\in\Gamma}g(\Pscr)=\B^3$;
   \item $g(\Int\,\Pscr)\cap\Int\,\Pscr=\emptyset$ for any nontrivial
         element $g\in\Gamma$, here $\Int\,\Pscr$ is the interior of
         $\Pscr$ in $\B^3$;
   \item For each side $S$ of $\Pscr$, there is another side $S'$ and an
         element $g\in\Gamma$ such that $g(S)=S'$;
   \item For any compact subset $K\subset\B^3$,
         $\{g\in\Gamma\ |\ g(\Pscr)\cap{}K\ne\emptyset\}$ is a finite set.
\end{enumerate}
\end{definition}
The relatively closed sets $Q_{\pm}=\overline{\Pscr}\cap\Omega_{\pm}$ are {\it fundamental domains}
of $\Gamma$ in $\Omega_{\pm}$, respectively. If $\Gamma$ is a {\qfg}, then we have
$\Pscr\cong{}Q_{+}\times\R$. Now we proceed with the following theorem which finds a closed {\ms} in
some submanifold of $M^3 =\H^3/\Gamma$.

\bt\label{key}
Let $\Lambda_\Gamma$ be the limit set of some torsion free {\qfg} $\Gamma$ such that $\B^3/\Gamma$ is
homotopic to a closed surface of genus at least two. If $\Pscr\subset\B^3$ is the fundamental polyhedra of the group
$\Gamma$ and $Q_{\pm}=\overline{\Pscr}\cap{}S_{\infty}^{2}\subset\Omega_{\pm}$ are the fundamental domains of
$\Gamma$ in $\Omega_{\pm}$, respectively, where $\Omega_{\pm}=S_{\infty}^{2}\setminus\Lambda_\Gamma$.
Suppose that $\Pi_{1},\ldots,\Pi_{n}$ are finite disjoint least area {\mct}s in $\overline{\B^{3}}$ that
satisfy the following conditions:
\begin{enumerate}
   \item The asymptotic boundary of each {\mct} must either both be contained in $Q_{+}$ or both be contained in $Q_{-}$;
   \item If $n=1$, then $\partial_{\infty}\Pi_{1}$ and $\Lambda_\Gamma$ are in good position; if $n\geq{}2$, then the {\abs}
   $\partial_{\infty}\Pi_{1},\ldots,\partial_{\infty}\Pi_{n}$ are in good position;
   \item $\Lambda_\Gamma$ is null-homotopic in the component of $\B^{3}\setminus\cup\Pi_{l}$ which contains
   $\Lambda_\Gamma$.
\end{enumerate}
Then there exists a $\Gamma$-invariant minimal disk with {\ab} $\Lambda_\Gamma$, which is disjoint from
$g(\Pi_{l})$, for all $g\in\Gamma$ and all $l\in\{1,\ldots,n\}$. Furthermore, the disk is embedded in $\B^3$.
\et

\bp
Our strategy is to remove some solid catenoids from $\B^3$ to obtain a submanifold of the {\qfm} $M^3 =\B^3/\Gamma$.

We let $\T_l$ be the solid catenoid bounded by $\Pi_l$ for each $l=1,\ldots,n$. By the assumptions and Lemma
~\ref{disjoint}, we know that these solid catenoids lie inside the {\fp} of the {\qfg} $\Gamma$, namely, $\T_{l}\subset\Pscr$
for each $l=1,\ldots,n$. Therefore, in what follows we abuse our notion to use $\T_{l}$
to denote the corresponding solid catenoid in $M^3$ as well.

From assumption of the theorem, we may assume that $\partial_{\infty}\Pi_{l}\subset{}Q_{+}$ for $l=1,\ldots,m$, and
$\partial_{\infty}\Pi_{l}\subset{}Q_{-}$ for $l=m+1,\ldots,n$, where $m\leq{}n$ are positive integers.

Since $\T_{l}\subset\Pscr$, we also have $g(\T_{l})\subset{}g(\Pscr)$ for all $g\in\Gamma$ and $l=1,\ldots,n$. Let
\be
   \X_{\infty}=\left\{\overline{\B^{3}\setminus
         \bigcup_{g\in\Gamma}\bigcup_{l=1}^{n}g(\T_l)}\right\}
         \setminus{}S_{\infty}^{2}\ ,
\ene
then $\X_{\infty}$ is a $\Gamma$-invariant subset of $\H^3$. Let $M_\infty=\X_{\infty}/\Gamma$, since
$\T_{l}\subset\Pscr$ for all $l=1,\ldots,n$, it is clear that $M_\infty$ is a submanifold of the {\qfm}
$M^3=\H^3/\Gamma$. By the construction, $M_\infty$ is a homogeneously regular $3$-manifold in the sense of
Morrey (\cite{Mor48}), whose boundary is mean convex {\wrt} the inward normal vector.

Now let $Y=\CC_{\Gamma}\cap{}M_\infty$, where $\CC_{\Gamma}$ is the {\cc} of the {\qfm} $M^3$, then $Y$ is a
compact {\htm} whose boundary is mean convex with respect to the inward normal vector (for the precise meaning of
positive {\mc} {\wrt} inward normal vector, see \S 7 of \cite{FHS83}).

By the Lemma ~\ref{Y} which we will prove later, this {\htm} $Y$ contains a closed surface that is incompressible. As
discussed in the subsection \S 2.2, we can now apply the results of Schoen-Yau, Sacks-Uhlenbeck, Freedman-Hass-Scott,
to conclude $Y$ also contains an embedded least area surface, denoted by $\Sigma$,  which is incompressible. Lifting this {\ms} $\Sigma$ to
$\H^3$, we obtain an embedded $\Gamma$-invariant minimal disk $\widetilde\Sigma$ in $\X_{\infty}$
with {\ab} $\Lambda_\Gamma$.
\ep

\begin{rem}In general, the universal cover, $\widetilde\Sigma$,
of an area minimizing surface $\Sigma$ in $M^3 = \H^3/\Gamma$
is not necessarily a least area minimal disk in $\H^3$ whose {\ab} is $\Lambda_\Gamma$.
\end{rem}
To complete the proof for the Theorem ~\ref{key}, we prove the following lemma:
\bl\label{Y}
There exists a closed incompressible surface $F$ contained in $Y$.
\el
\bp
The {\cc} of a {\qfm} has two boundary components, namely we can write
$\partial\CC_{\Gamma}=F_{+}\cup{}F_{-}$, where $F_{\pm}$ is the pleated surface that faces the surface
$\Omega_{\pm}/\Gamma$, respectively.

Let $W_{\pm}$ be the submanifold of $M^3$ that is bounded by $F_{\mp}$ and the conformal infinity
$\Omega_{\pm}/\Gamma$, respectively. Obviously $\CC_{\Gamma}=W_{+}\cap{}W_{-}$.
Then $\partial{}W_{\pm}=F_{\mp}$ is convex {\wrt} the inward normal vector.
For all $l=1,\ldots,n$, since the catenoid
$\Pi_{l}=\partial\T_{l}$ is a least area {\ms},
by \cite[Theorem 1]{MY82a}, $\Pi_{l}\subset{}W_{+}$ if $1\leq{}l\leq{}m$, and
$\Pi_{l}\subset{}W_{-}$ if $m+1\leq{}l\leq{}n$.
Therefore $\T_{l}\subset{}W_{+}$ if $1\leq{}l\leq{}m$, and
$\T_{l}\subset{}W_{-}$ if $m+1\leq{}l\leq{}n$.

Since $\Lambda_\Gamma$ is null-homotopic in the component of
$\B^{3}\setminus\cup\Pi_{l}$ which contains $\Lambda_\Gamma$, there
exists a disk $D$ that is asymptotic to $\Lambda_\Gamma$ and that is
disjoint from $\{\T_1,\ldots,\T_n\}$. In particular, $D\cap\Pscr$ is disjoint
from $\{\T_1,\ldots,\T_n\}$. Then we can perturb $D\cap\Pscr$ to
get a disk $\Delta\subset\Pscr$ so that $\Delta$ is still disjoint
from $\{\T_1,\ldots,\T_n\}$ and
\beq
   \Dcal=\bigcup_{g\in\Gamma}g(\Delta)
\eeq
is a $\Gamma$-invariant disk asymptotic to $\Lambda_\Gamma$, and it is disjoint from $g(\T_l)$ for
all $g\in\Gamma$ and all $l\in\{1,\ldots,n\}$.

Let $F'=\Dcal/\Gamma$, then clearly $F'$ is a closed incompressible surface. If $F'$ is contained in $Y$, then we are done.
Otherwise, by the above discussion, $F_{+}$ is disjoint from $\T_{l}$ for all $m+1\leq{}l\leq{}n$ and $F_{-}$ is disjoint from
$\T_{l}$ for all $1\leq{}l\leq{}m$, so $F'$ is isotopic to an incompressible surface $F\subset{}Y$ via an isotopy
$\varphi_{t}:M_{\infty}\to{}M_{\infty}$ such that $\varphi_{0}=\id_{M_\infty}$ and each $\varphi_{t}$ is a
homeomorphism, where $0\leq{}t\leq{}1$.
\ep
\subsection{A topological lemma}

We will prove a topological lemma in this subsection (see Lemma \ref{lem:essential Jordan curve}). This will be very useful 
when we show the distinctiveness of {\ms}s from the construction.

Let us recall some notations. Let $\Lambda$ be a {\jc} on the sphere at infinity $S_\infty^2$, and 
$\Omega_{\pm}=S_\infty^2\setminus\Lambda$. Suppose there are two disjoint solid {\mct}s $\T_{\pm}$ with boundary circles 
$\partial_\infty\T_{\pm}=\Delta'_{\pm}\cup\Delta_{\pm}''\subset\Omega_{\pm}$. Let $\delta_{\pm}\subset\B^3\cup{}S_{\infty}^2$ 
be two (unknotted) simple closed curves constructed as follows: take two pints $p'_{+}\in\Delta'_{+}$ and $p''_{+}\in\Delta''_{+}$. Let
$\alpha_{+}\subset\T_{+}$ be a curve with asymptotic endpoints $p'_{+}$ and $p''_{+}$, and let $\beta_{+}\subset\Omega_{+}$ 
be a simple curve with endpoints $p'_{+}$ and $p''_{+}$. Let $\delta_{+}=\alpha_{+}\cup\beta_{+}$. Similarly we can define the loop 
$\delta_{-}$.

\begin{figure}[htbp]
\centering
\begin{minipage}[t]{0.45\linewidth}
\centering
\includegraphics[scale=0.5]{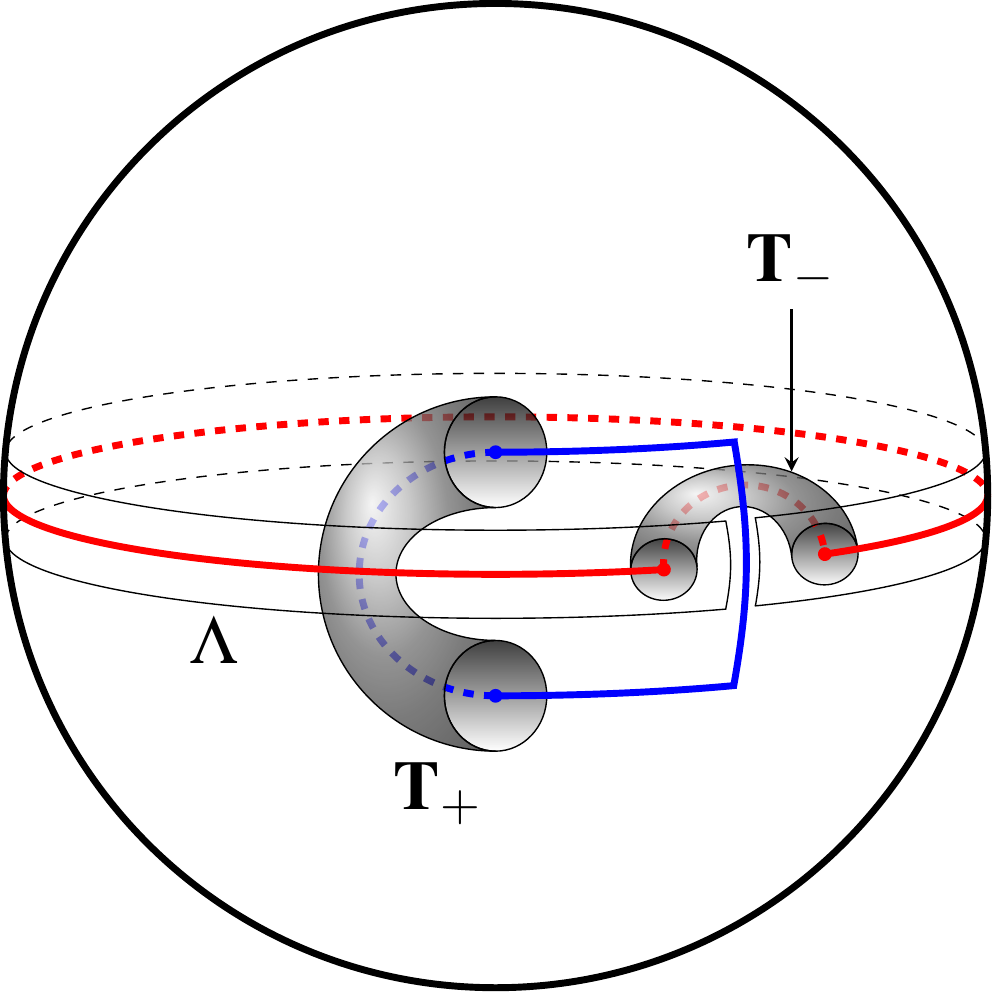}
\end{minipage}%
\begin{minipage}[t]{0.45\linewidth}
\centering
\includegraphics[scale=0.5]{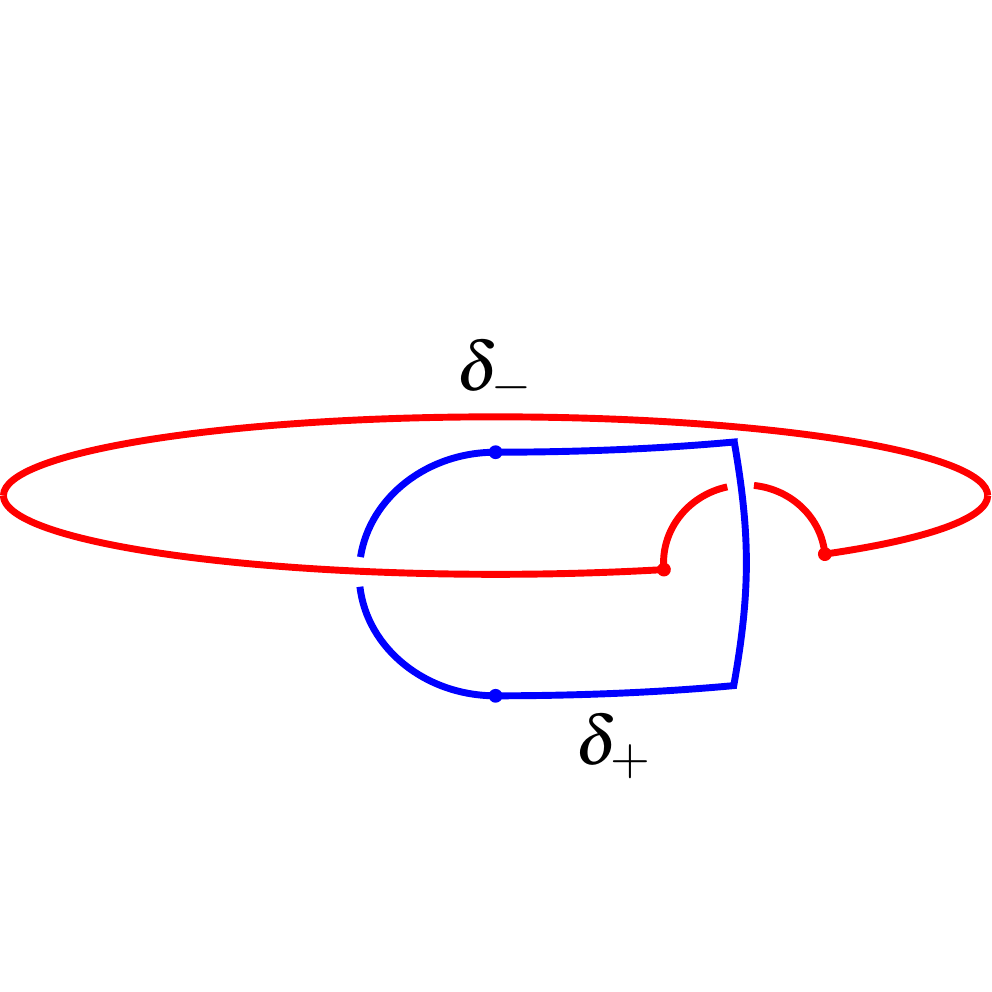}
\end{minipage}
\caption{Linked loops in $\B^3\cup{}S_{\infty}^2$}
\label{fig:essential Jordan curve}
\end{figure}

We say that solid {\mct}s $\T_{+}$ and $\T_{-}$ are {\em linked} (or {\em unlinked}) in $\B^3\cup{}S_{\infty}^2$ if $\delta_{+}$ and 
$\delta_{-}$ are linked (or unlinked) loops. See Figure \ref{fig:essential Jordan curve} for the linked case.

\begin{lem}\label{lem:essential Jordan curve}
We have the following statements:
\begin{enumerate}
   \item If solid {\mct}s $\T_{+}$ and $\T_{-}$ are linked in $\B^3\cup{}S_{\infty}^2$,
         then the Jordan curve $\Lambda$ is essential in the space
         $\overline{\B^3\setminus(\T_{+}\cup\T_{-})}$.
   \item If $\T_{+}$ and $\T_{-}$ are unlinked in $\B^3\cup{}S_{\infty}^2$,
         then $\Lambda$ is null-homotopic in the space
         $\overline{\B^3\setminus(\T_{+}\cup\T_{-})}$.
         
\end{enumerate}
\end{lem}

\begin{figure}[htbp]
\centering
\begin{minipage}[t]{0.5\linewidth}
\centering
\includegraphics[scale=0.85]{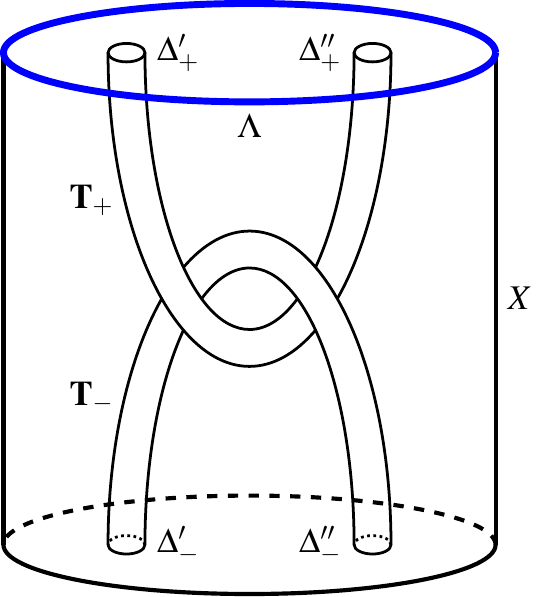}
\end{minipage}%
\begin{minipage}[t]{0.5\linewidth}
\centering
\includegraphics[scale=1]{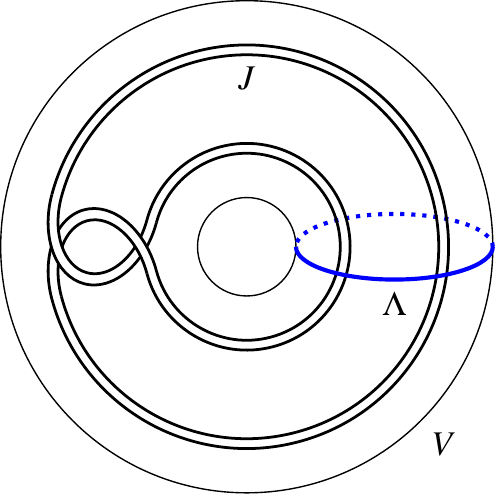}
\end{minipage}%
\caption{Jordan curve $\Lambda$ is essential in
$\overline{\B^3\setminus(\T_{+}\cup\T_{-})}$}\label{fig:essential Jordan curve2}
\end{figure}

\begin{rem}If the boundary circles $\partial_{\infty}\T_{+}$ and $\partial_{\infty}\T_{-}$ are contained
in the same component of $S_\infty^2\setminus\Lambda$, then $\Lambda$ is always
null-homotopic in $\overline{\B^3\setminus(\T_{+}\cup\T_{-})}$.
\end{rem}
\bp(of Lemma \ref{lem:essential Jordan curve})
We will only prove the statement (1) as the other follows easily. We consider the equivalent figure of
Figure \ref{fig:essential Jordan curve} (see Exercise I.6 in \cite[Chapter 3]{Rol90}): Let $X$ denote the solid 
cylinder in Figure \ref{fig:essential Jordan curve2}, and let $V$ be the quotient space of $X$ obtained by 
gluing the top and the bottom disks of $X$. We also identify $\Delta'_{+}$ with $\Delta'_{-}$ and $\Delta''_{+}$ 
with $\Delta''_{-}$ in the quotient space $V$, then $J=\T_{+}\cup\T_{-}$ is a solid torus contained in $V$.

If the {\jc} $\Lambda$ is null-homotopic in $X\setminus(\T_{+}\cup\T_{-})$, then it is also null-homotopic in 
$V\setminus{}J$. But this is impossible according to Proposition G.3 in \cite[Chapter 3]{Rol90}.
Therefore $\Lambda$ is not contractible in $\overline{\B^3\setminus(\T_{+}\cup\T_{-})}$.

\ep

\subsection{Proof of the main theorem} In this subsection, we use results of the previous subsections to
prove our main Theorem ~\ref{main}, namely, for a prescribed integer $N$, we construct a {\qfm}
$M^3 = \H^3/\Gamma$ such that it contains at least $2^N$ many closed incompressible {\ms}s.

All figures in this subsection correspond to the case $N = 3$ for simplicity.

\bp[\bf{Proof of Theorem \ref{main}}] The construction consists of four steps which we now describe.
\vskip 0.1in
\noindent
\textbf{Step 1}:  Consider the unit ball $\B^3$ in the $xyz$-space. We divide the interval $[-1,1]$ on the $z$-axis into
$N+1$ subintervals of equal length. Let $\varepsilon>0$ be a sufficiently small number and $H^{\pm}_{i}$ be the
(Euclidean) horizontal planes
\beq
   z=1-\frac{2i}{N+1}\pm\varepsilon,
   \quad{}1\leq{}i\leq{}N\ .
\eeq
Let $C^{\pm}_{i}=H^{\pm}_{i}\cap{}S_{\infty}^{2}$ be a pair of parallel
circles on $S_{\infty}^{2}$, $1\leq{}i\leq{}N$.

\begin{figure}[htbp]
\begin{minipage}[t]{0.75\linewidth}
\centering
\includegraphics[scale=0.5]{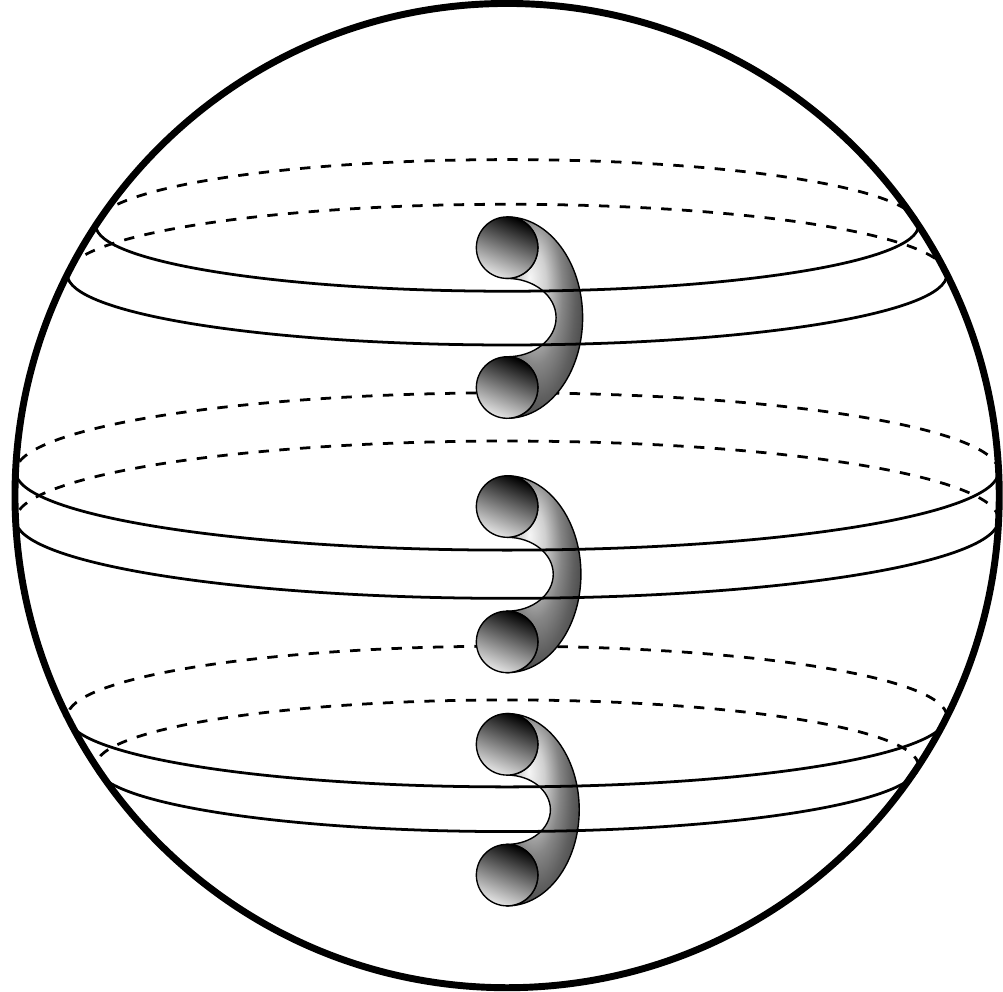}
\caption{$2N$ Parallel circles $(N=3)$}\label{fig:parallel circles}
\end{minipage}%
\end{figure}

By Lemma~\ref{0dist} and Theorem~\ref{Gomes}, we may assume that $\varepsilon$ is sufficiently small
so that we can find a pair of circles above $C^{+}_{i}$ and below $C^{-}_{i}$, respectively, that bound a spherical
{\mct} for each $1\leq{}i\leq{}N$ (see Figure~\ref{fig:parallel circles}). By the Theorem ~\ref{wang}, we may
assume each {\mct} is of least area.

\begin{figure}[htbp]
\begin{minipage}[t]{0.75\linewidth}
\centering
\includegraphics[scale=0.5]{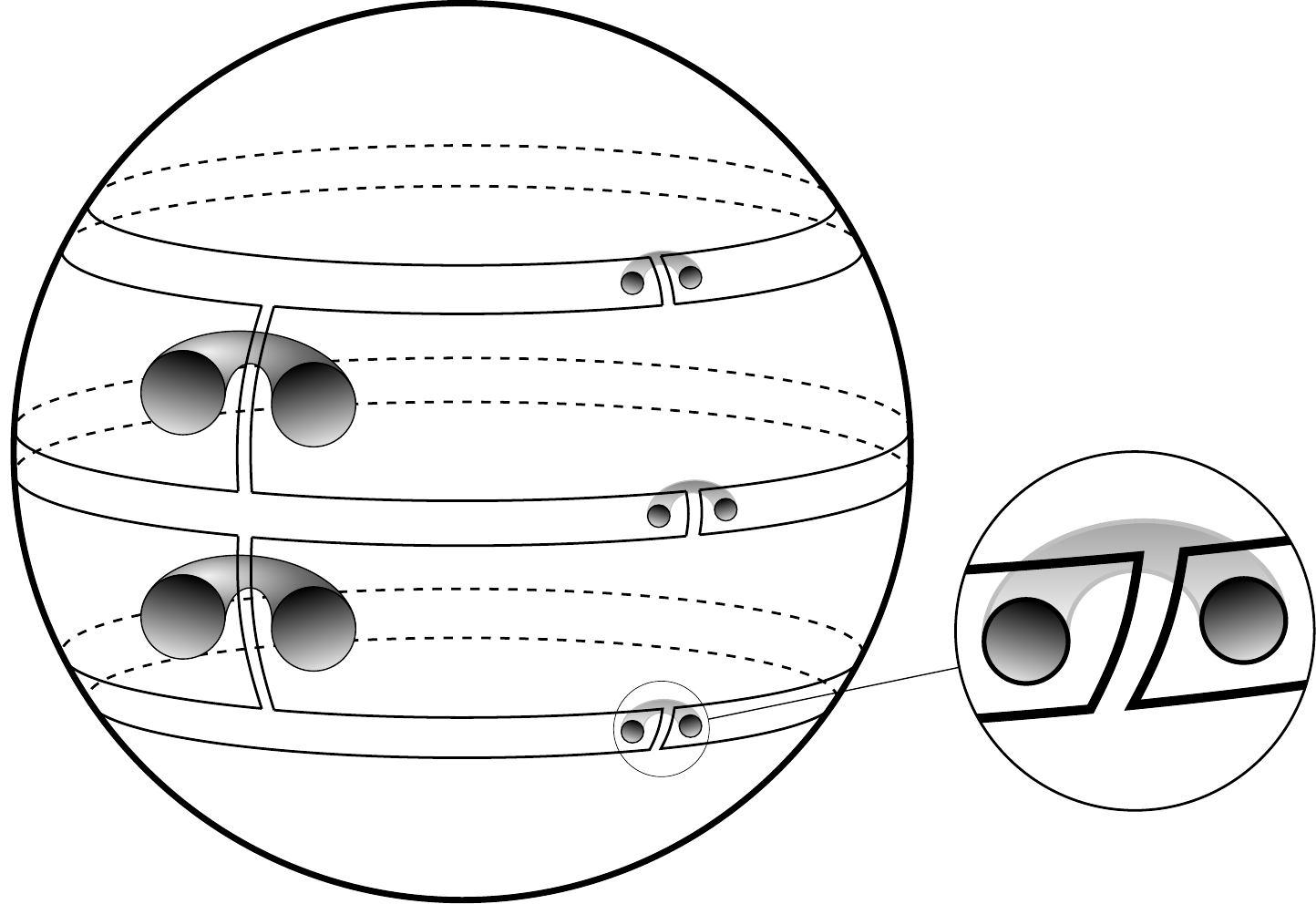}
\caption{$2N-1$ Narrow bridges $(N=3)$}\label{fig:narrow bridges}
\end{minipage}
\end{figure}

\vskip 0.1in
\noindent
\textbf{Step 2}: We connect each pair of parallel circles $C^{\pm}_{i}$ by a narrow bridge $B_i$, ($1\leq{}i\leq{}N$),
and connect each pair of parallel circles $C^{-}_{i}$ and $C^{+}_{i+1}$ by a narrow bridge $B_{j}'$
($1\leq{}j\leq{}N-1$), then we obtain a piecewise smooth {\jc} $\Lambda\subset{}S_{\infty}^{2}$
(see Figure~\ref{fig:narrow bridges}).

By Lemma~\ref{0dist} and Theorem ~\ref{wang}, we may assume that the $2N-1$ bridges are again
sufficiently narrow so that we can find a pair of circles around each bridge that bound a spherical
least area {\mct}.

\begin{figure}[htbp]
\begin{minipage}[t]{0.75\linewidth}
\centering
\includegraphics[scale=0.5]{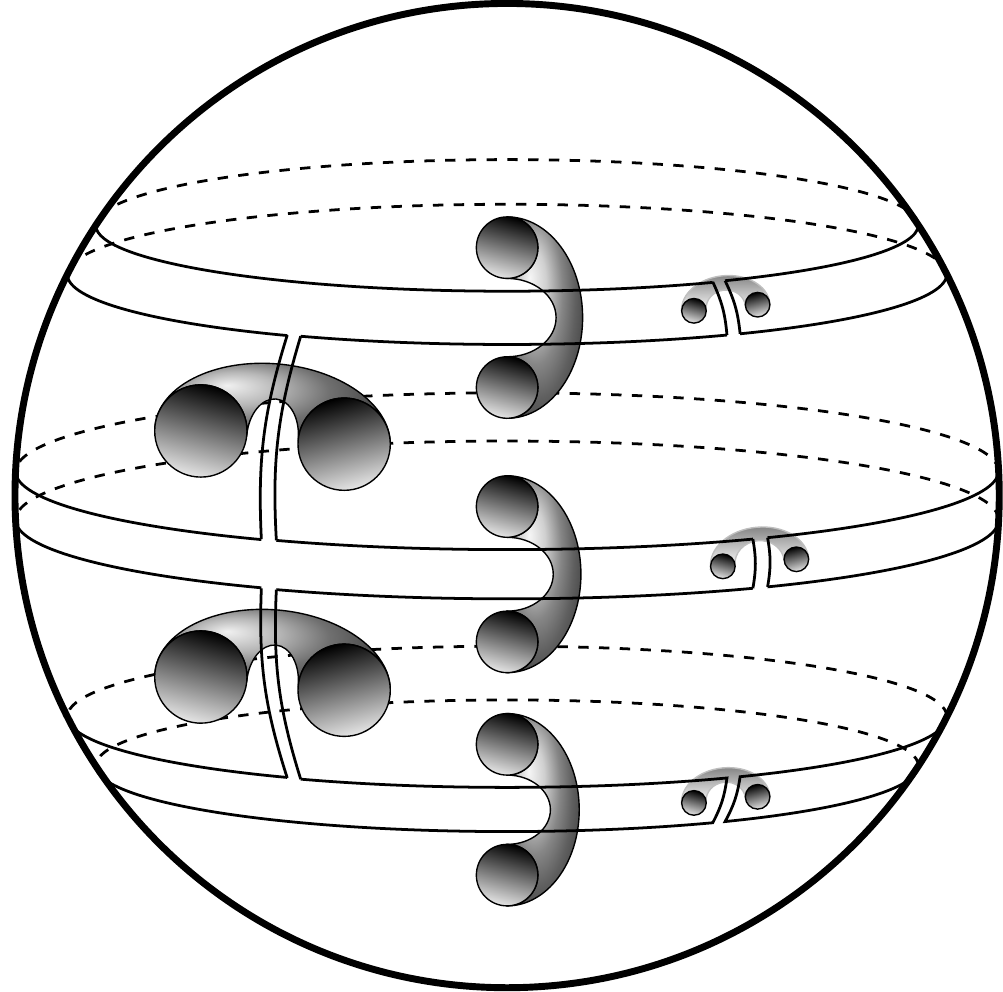}
\caption{$3N-1$ Catenoids $(N=3)$}\label{fig:catenoids}
\end{minipage}
\end{figure}

By the construction, there are total of $6N-2$ disjoint circles, and they are in good position. We also obtain
$3N-1$ solid {\mct}s in $\B^3$, for each pair of disjoint circles next to the same bridge.
\vskip 0.1in
\noindent
\textbf{Step 3}: Now we construct the {\qfg} $\Gamma$ by the method in subsection \S 3.2. In particular, we
cover the Jordan curve $\Lambda$ in Step 2 by $2L$ small disks $\{\Delta_1,\Delta_2,\ldots,\Delta_{2L}\}$
on $S_{\infty}^{2}$. We can choose these disks sufficiently small such that $\cup\Delta_{i}$ is disjoint from
the $6N-2$ circles in step two (see Figure ~\ref{fig:catenoids}).

Now we are in position to apply Theorem ~\ref{inversion} to use the inversions {\wrt} the circles
$\{\partial\Delta_1,\partial\Delta_2,\ldots,\partial\Delta_{2L}\}$ to construct a torsion free {\qfg} $\Gamma$
whose limit set $\Lambda_{\Gamma}$ is contained in $\cup_{i=1}^{2L}\Delta_{i}$.

It is easy to see that $\cup_{i=1}^{2L}\Delta_{i}$ separates the sphere $S_{\infty}^{2}$ into two
components, denoted by $E_{\pm}$. By the construction, it is clear that for each pair of the $6N-2$ circles,
either both are contained in $E_{+}$ or both are in $E_{-}$. Therefore for each pair of these circles, either both
are in $Q_{+}$ or both are in $Q_{-}$, where $E_{\pm}\subset{}Q_{\pm}\subset\Omega_{\pm}$ are the {\ab} of a
fundamental polyhedron $\Pscr$ of $\Gamma$.

\begin{figure}[htbp]
\centering
\begin{minipage}[t]{0.45\linewidth}
\centering
\includegraphics[scale=0.425]{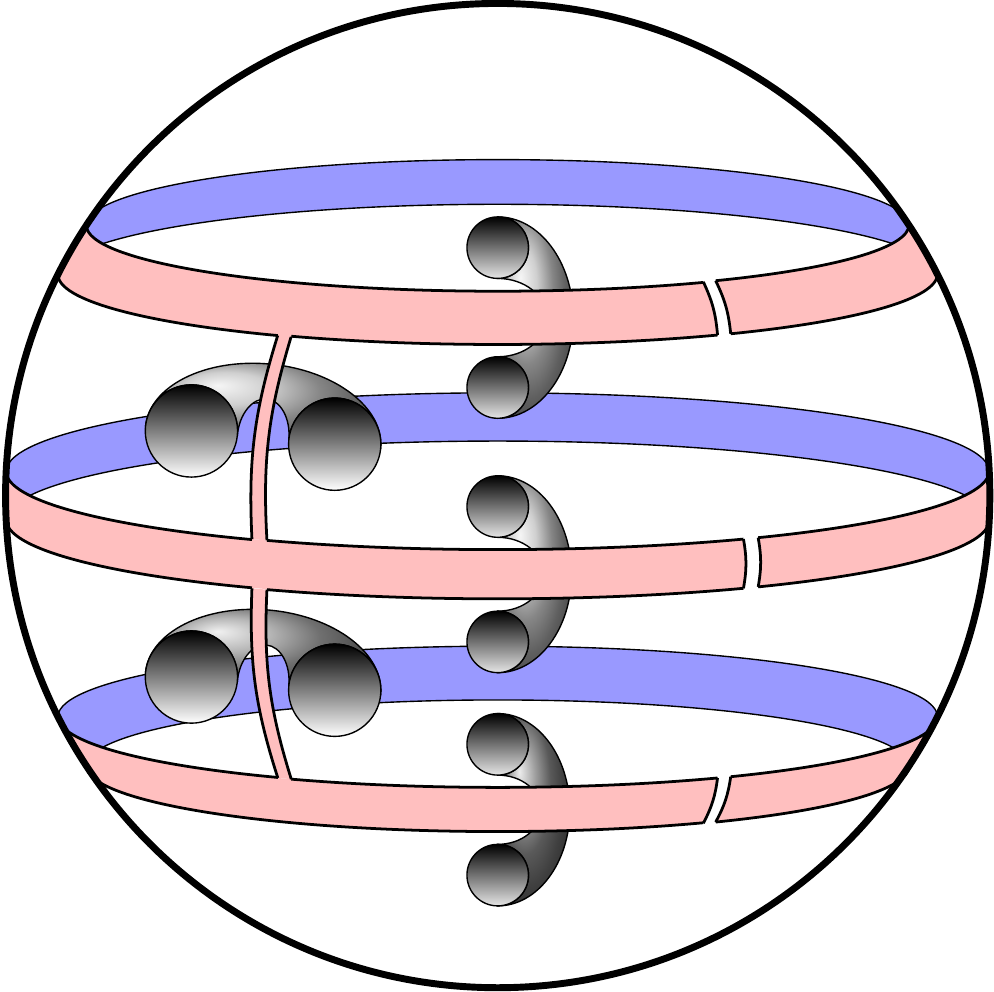}
\end{minipage}%
\begin{minipage}[t]{0.45\linewidth}
\centering
\includegraphics[scale=0.425]{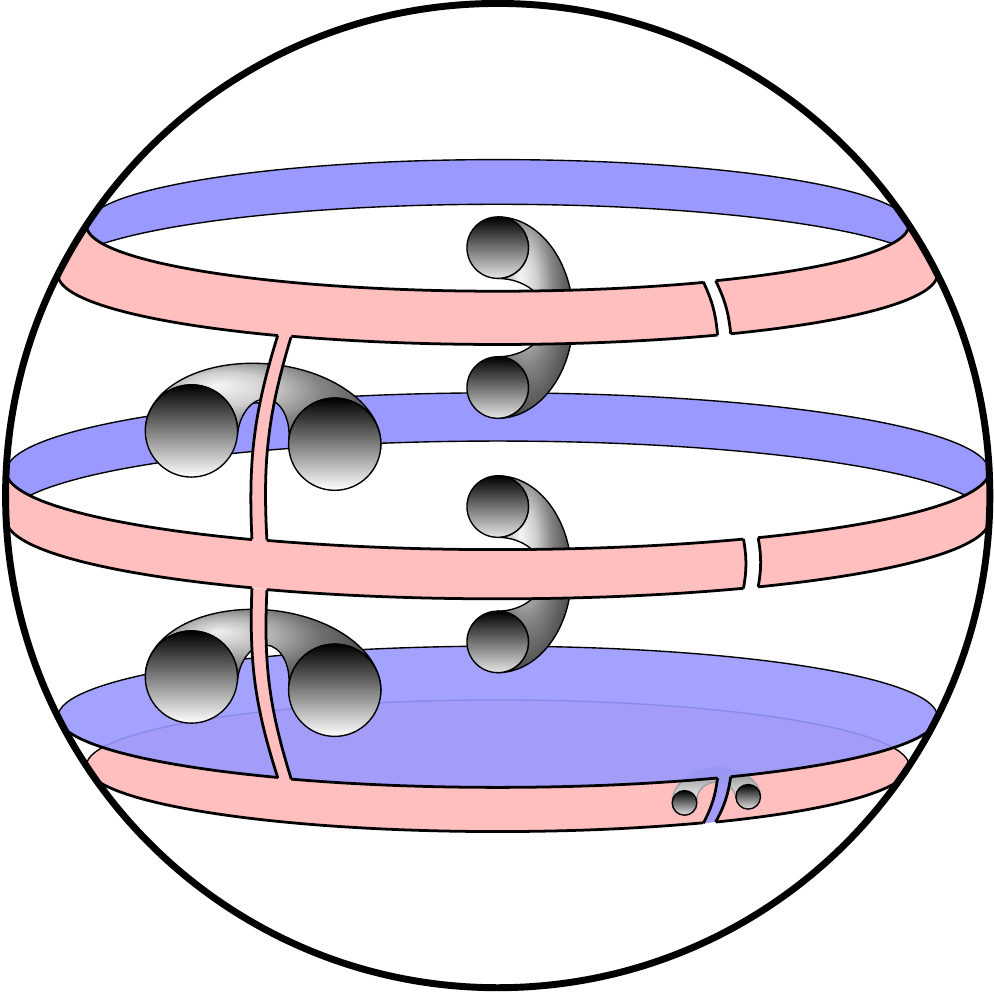}
\end{minipage}\\[20pt]%
\begin{minipage}[t]{0.45\linewidth}
\centering
\includegraphics[scale=0.425]{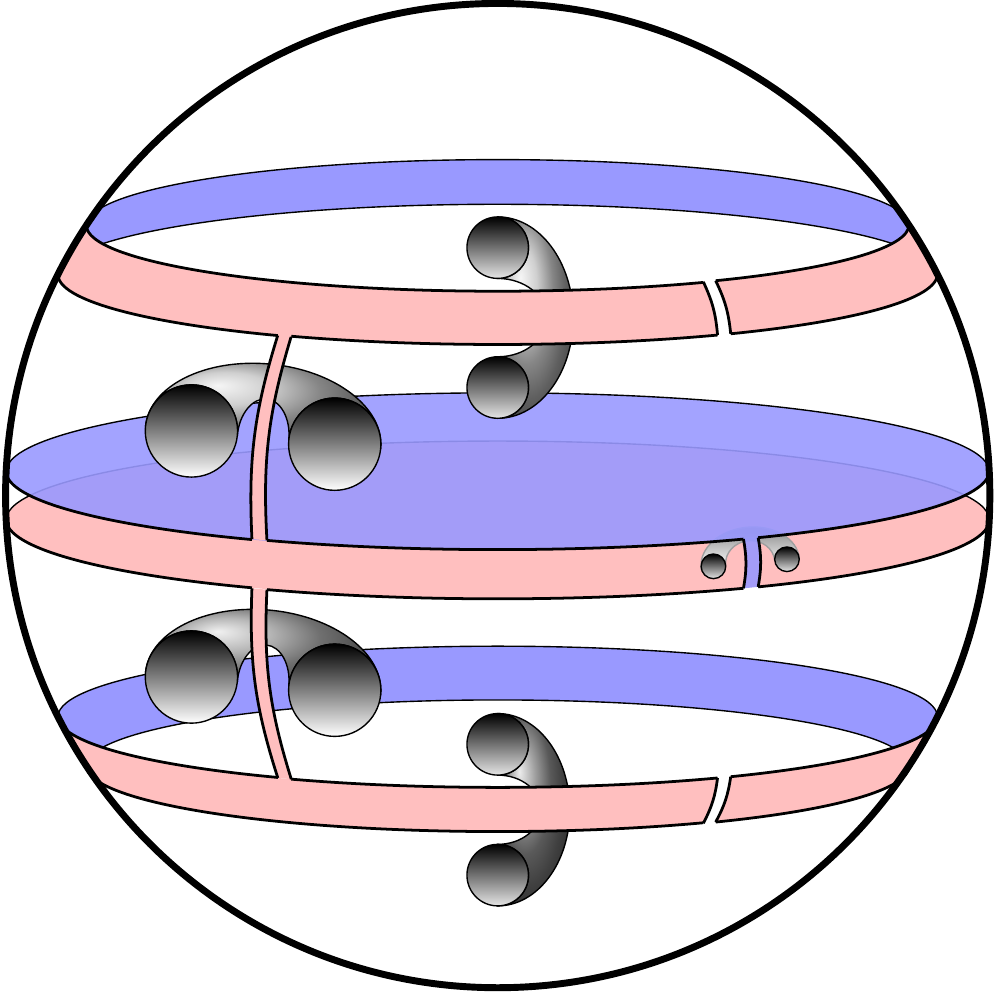}
\end{minipage}%
\begin{minipage}[t]{0.45\linewidth}
\centering
\includegraphics[scale=0.425]{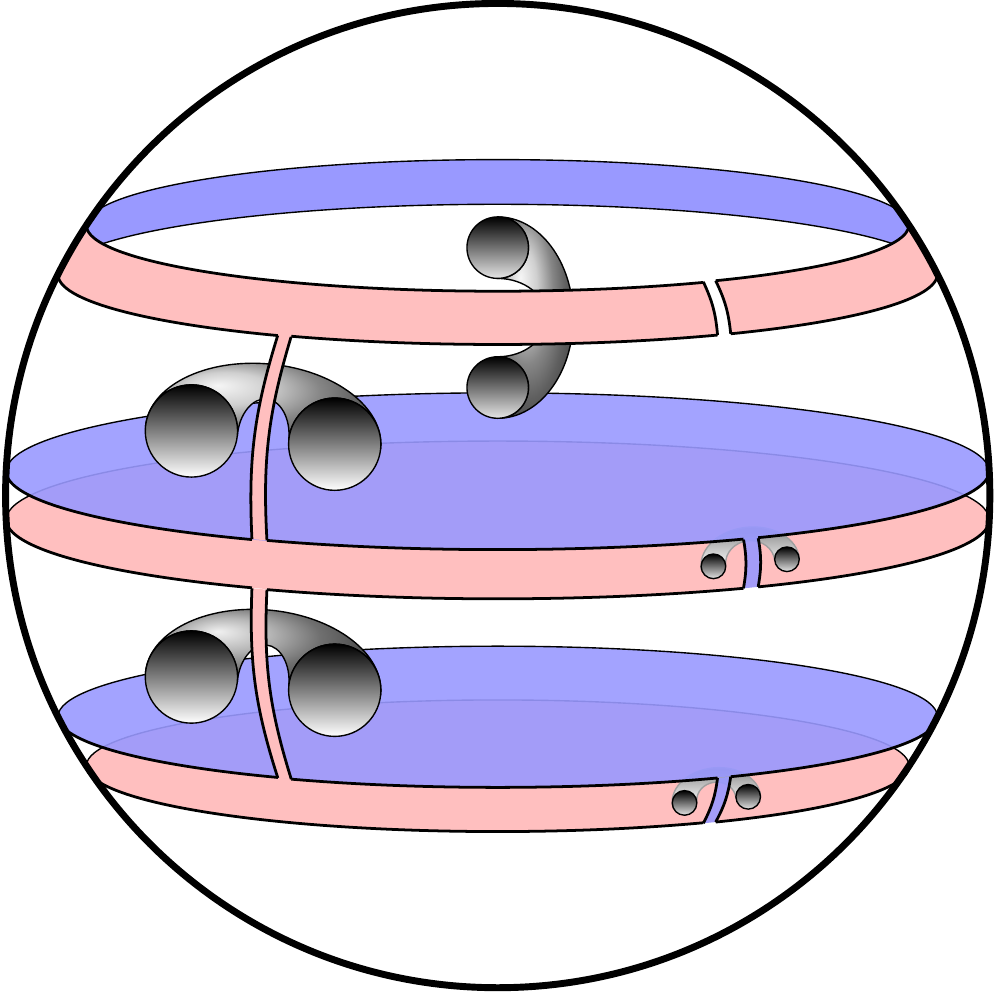}
\end{minipage}\\[20pt]%
\begin{minipage}[t]{0.45\linewidth}
\centering
\includegraphics[scale=0.425]{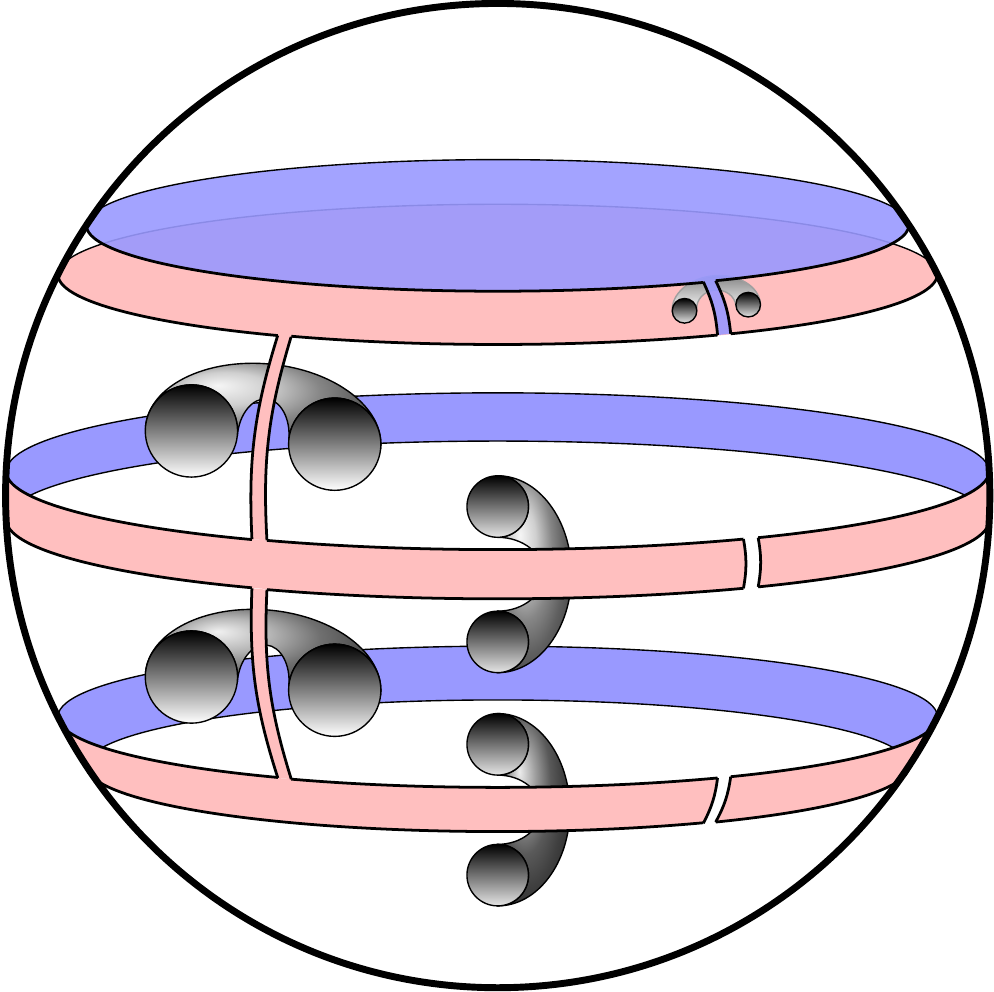}
\end{minipage}%
\begin{minipage}[t]{0.45\linewidth}
\centering
\includegraphics[scale=0.425]{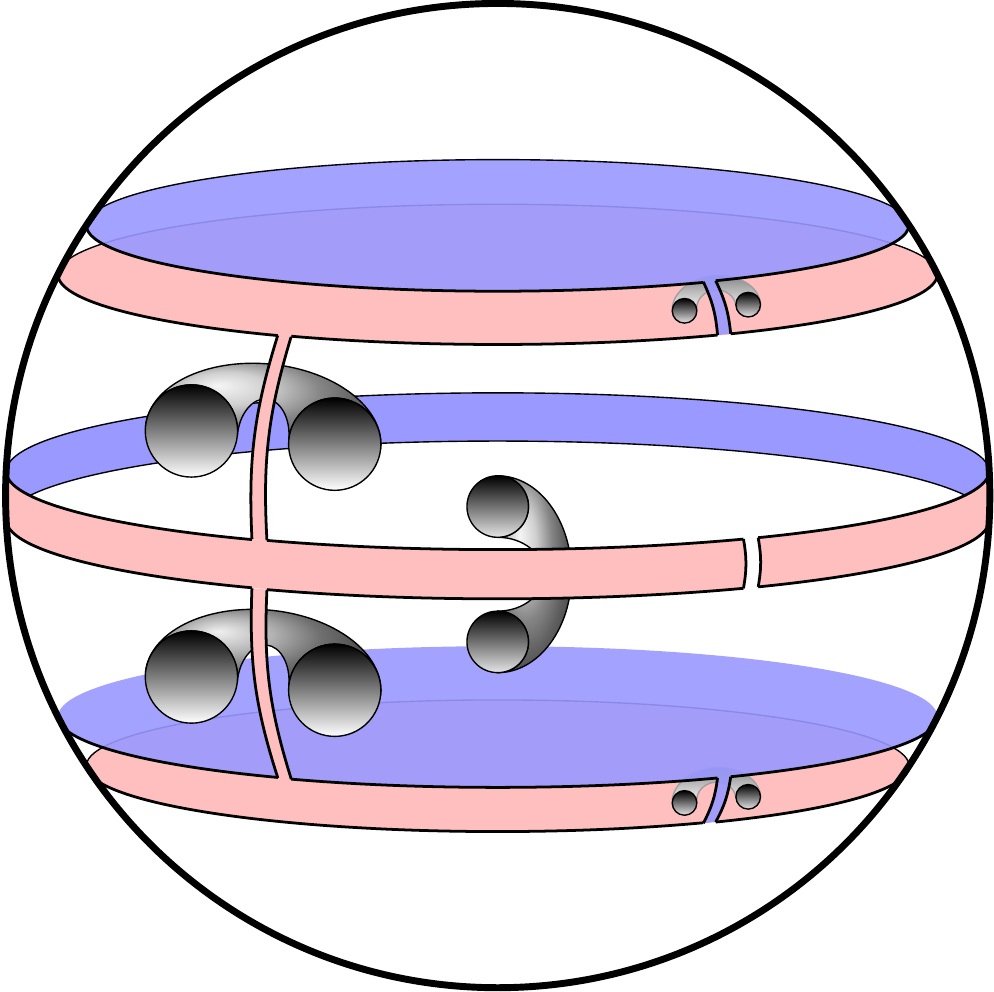}
\end{minipage}\\[20pt]%
\begin{minipage}[t]{0.45\linewidth}
\centering
\includegraphics[scale=0.425]{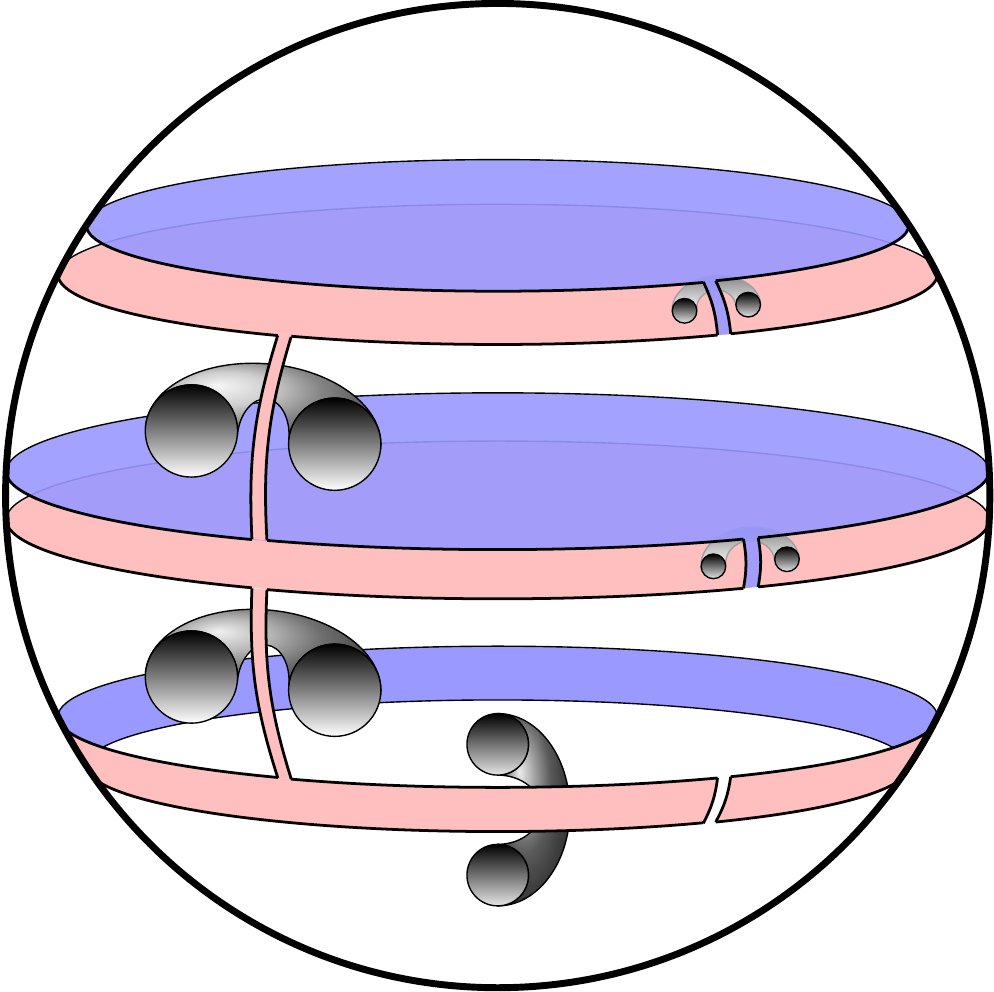}
\end{minipage}%
\begin{minipage}[t]{0.45\linewidth}
\centering
\includegraphics[scale=0.425]{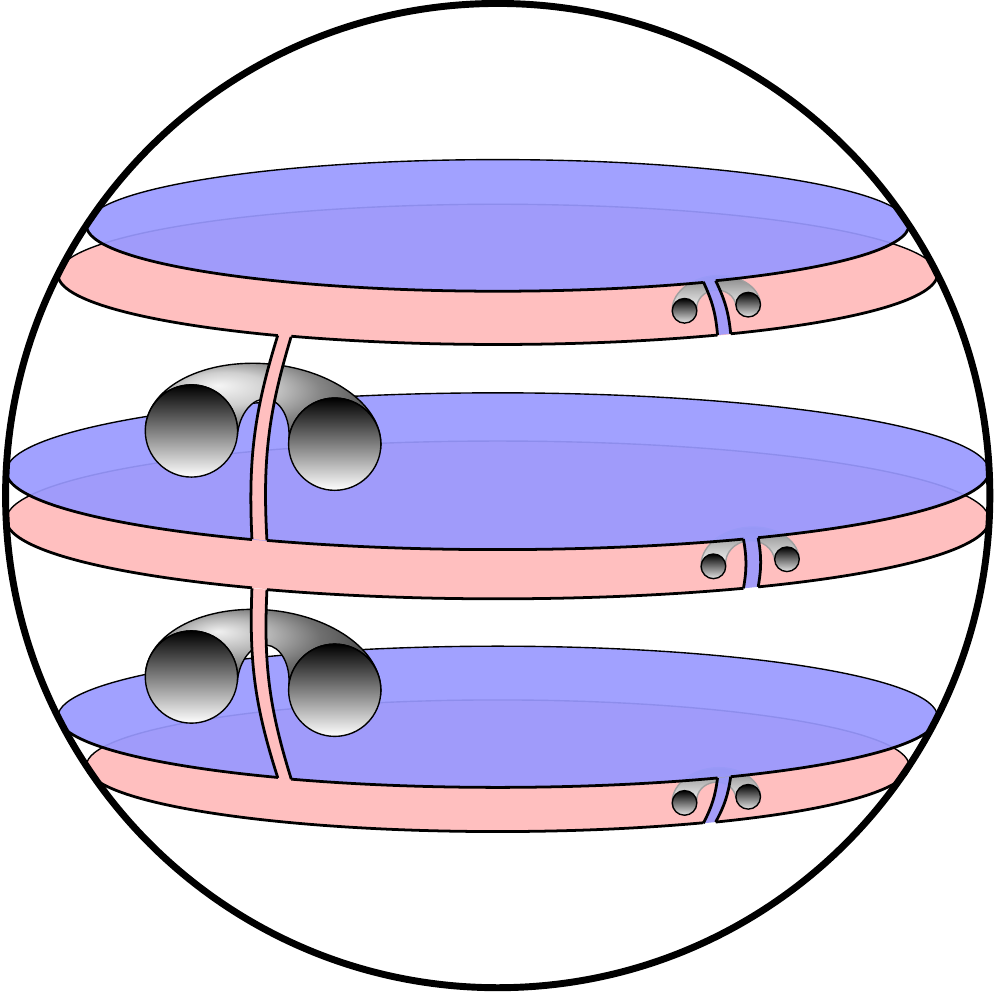}
\end{minipage}%
\caption{$2^N$ Minimal Surfaces ($N=3$)}\label{fig:all surfaces}
\end{figure}

\vskip 0.1in
\noindent
\textbf{Step 4}: Now we can choose the $3N-1$ pairs of circles in different ways to form $2N-1$ {\mct}s
(see Figure~\ref{fig:all surfaces} for the case when $N=3$). Let us fix some notations: each {\jc} in 
Figure~\ref{fig:all surfaces} represents the limit set $\Lambda_\Gamma$ of the {\qf} group $\Gamma$.
Let the circles next to bridges $B'_j$ ($1\leq j \leq N-1$) be fixed.
Therefore the $N-1$ {\mct}s bounded by these circles are fixed as well. Now we are left with
$2N$ {\mct}s: $N$ bounded by each pair of {\tdc} next to the bridge $B_i$ ($1\leq i \le N$), and the
other $N$, each bounded by a circle above the circle $C_i^{+}$ and a circle below the circle $C_i^{-}$
($1\leq i \le N$), as in the Figure ~\ref{fig:parallel circles}. Namely, for each bridge $B_i$, $1\leq i \le N$,
there are two choices of {\mct}s. Hence we have a choice of $2^N$ many different ways to obtain $2N-1$
{\mct}s, which can be seen from the following: let $\Pi'_i, \Pi_i^B$, and $\Pi_i^C$ be the {\mct}s asymptotic
to the circles near the bridge $B'_i$, the bridge $B_i$ and the circle $C_i$, respectively. For each
function $\imath: \{1,\cdots,N\} \to \{B,C\}$, let $\Pi = \bigcup_{i=1}^N \Pi_i^{\imath(i)}$. We see that for each
$\imath$, $\Pi_{\imath}$ is a collection of $N$ {\mct}s. Such a collection, together with
$\bigcup_{i=1}^{N-1} \Pi'_i$, is then a collection of $2N-1$ {\mct}s. Since there are $2^N$ choices of $\imath$,
there are $2^N$ such collections. For example, in the Figure~\ref{fig:all surfaces} where $N=3$, there are $8$
total choices, each case leading to  $5$ {\mct}s.

\begin{figure}[htbp]
\begin{minipage}[t]{0.95\linewidth}
\centering
\includegraphics[scale=0.5]{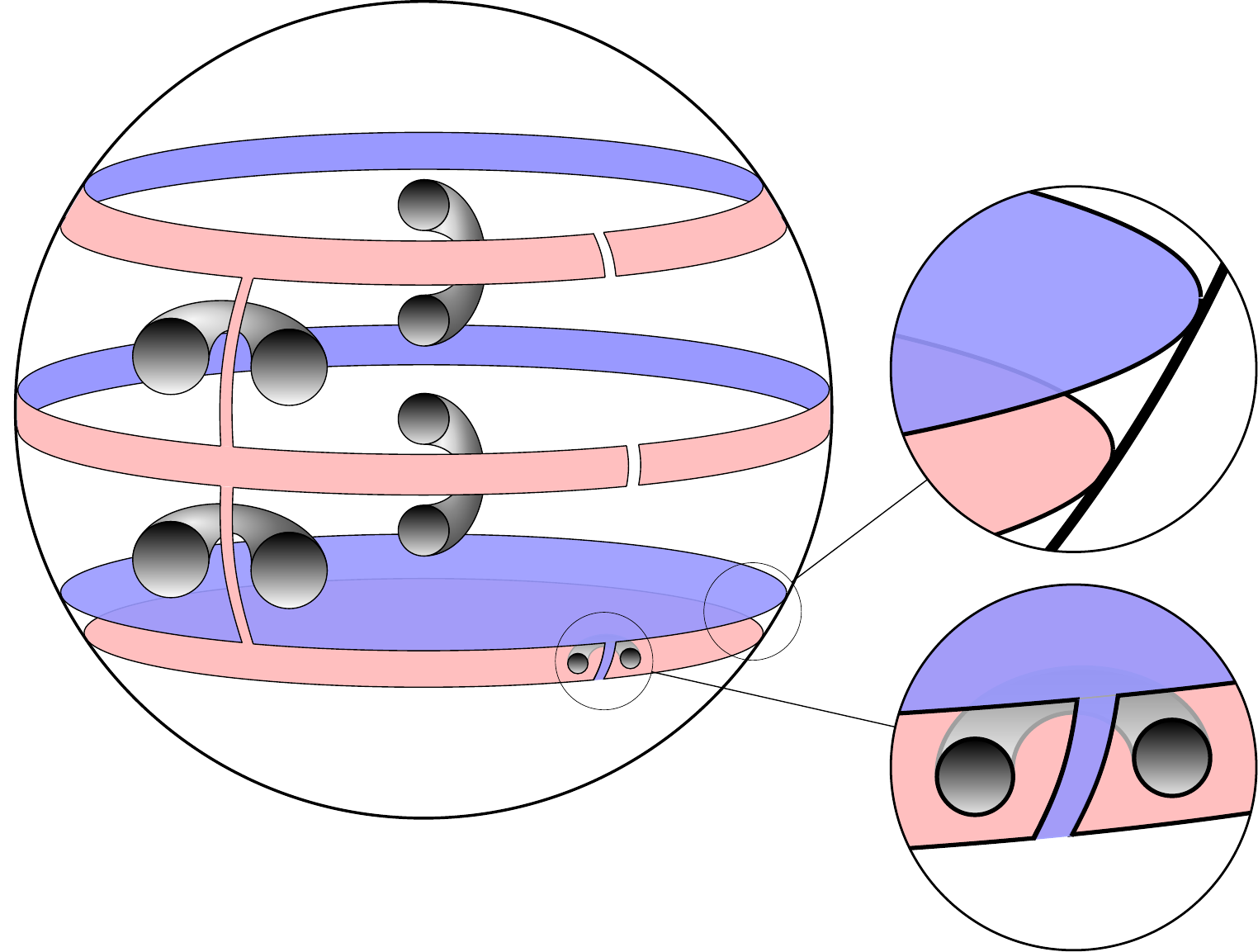}
\caption{$\Lambda_\Gamma$ is null-homotopic in the subspace of $\B^3$
with $2N-1$ solid catenoids removed $(N=3)$}\label{fig:null-homotopic}
\end{minipage}%
\end{figure}

Now for each case, it is easy to verify that $\Lambda_\Gamma$ bounds a disk, which consists of
strips and/or disks connected by narrow bridges, and are disjoint from $2N-1$ solid catenoids, therefore it is 
null-homotopic in the subspace of $\B^3$ with $2N-1$ solid catenoids removed. For example, when $N=3$, we 
consider the top right figure in Figure \ref{fig:all surfaces}, it is easy to see that $\Lambda_\Gamma$ bounds a disk
consisting of two strips and two disks connected by $3$ narrow bridges, therefore $\Lambda_\Gamma$ is 
null-homotopic in the subspace of $\B^3$ with $5$ solid catenoids removed (see Figure \ref{fig:null-homotopic}).

\vskip 0.1in

Now we apply Theorem ~\ref{key}, for each case, there is a $\Gamma$-invariant minimal disk with {\ab} 
$\Lambda_\Gamma$. Therefore we obtain $2^N$ $\Gamma$-invariant minimal disks in $\B^3$, denoted by 
$\Dcal_{1},\ldots,\Dcal_{2^{N}}$, which all share the {\ab} $\Lambda_\Gamma$.

We claim that these minimal disks $\Dcal_{1},\ldots,\Dcal_{2^{N}}$ are distinct. To illustrate this, 
we once again 
consider the case when $N=3$. We count the {\ms}s in Figure \ref{fig:all surfaces} from left to right and from 
top to bottom. Recall that each minimal disk $\Dcal_{i}$ ($i=1,\ldots,8$) is isotopic to the disk in the $i$-th subspace 
of $\B^3$ with $5$ catenoids removed in Figure \ref{fig:all surfaces}. We can now show, for example, that $\Dcal_2$ 
(see Figure \ref{fig:null-homotopic}) is distinct from other seven minimal disks. Indeed, the disk $\Dcal_2$ is distinct 
from $\Dcal_1$, since $\Dcal_1$ intersects the bottom catenoid in Figure \ref{fig:null-homotopic} by 
Lemma \ref{lem:essential Jordan curve}; similarly, by Lemma \ref{lem:essential Jordan curve}, each 
$\Dcal_i$, $i=3,\ldots,8$, intersects either the top catenoid, the middle catenoid or both in Figure \ref{fig:null-homotopic}. 
Therefore $\Dcal_{1},\ldots,\Dcal_{8}$ are distinct minimal disks. The general case is similar,  therefore these minimal 
disks $\Dcal_{1},\ldots,\Dcal_{2^{N}}$ from our construction are distinct.

Now we set $\Sigma_{i}=\Dcal_{i}/\Gamma$, $i=1,\ldots,2^N$, then $\Sigma_{1},\ldots,\Sigma_{2^{N}}$ are 
distinct embedded incompressible {\ms}s in the {\qfm} $M^3 = \B^3/\Gamma$. This completes the proof of our main theorem.
\ep

In Step 4, the minimal disks $\Dcal_{1},\ldots,\Dcal_{2^{N}}$ share the {\ab} $\Gamma$, which is the
{\jc} justifying Theorem 1.2.

\begin{rem}[On the genus of the surface]
By our construction, we use $2L$ small disks of injective radii $<\delta$ to cover the {\jc} $\Lambda$, then by the
Theorem ~\ref{inversion}, we can construct a torsion free {\qfg} $\Gamma$ so that
$\B^3/\Gamma\cong{}F\times\R$, where $F$ is a closed surface with $\genus(F)=L-1$. One sees that
$\Length(\Lambda)/\delta$ is large, and therefore so is the genus of $F$, and $L\to\infty$ as $N\to\infty$.
\end{rem}


\bibliographystyle{amsalpha}
\bibliography{ref_count}
\end{document}